\documentclass[11pt]{article}
\usepackage[psamsfonts]{amsfonts}
\usepackage{amsmath,amssymb,amsthm}

\usepackage{graphicx}
\usepackage{eepic}

\newcommand{\bt}   {\begin{theorem}}
\newcommand{\et}   {\end  {theorem}}
\newcommand{\bl}   {\begin{lemma}}
\newcommand{\el}   {\end  {lemma}}
\newcommand{\bp}   {\begin{prop}}
\newcommand{\ep}   {\end  {prop}}
\newcommand{\bc}   {\begin{cor}}
\newcommand{\ec}   {\end  {cor}}
\newcommand{\bd}   {\begin{defn}}
\newcommand{\ed}   {\end  {defn}}

\newcommand{\ba}   {\begin{array}}
\newcommand{\ea}   {\end  {array}}
\newcommand{\be}   {\begin{enumerate}}
\newcommand{\ee}   {\end  {enumerate}}
\newcommand{\bi}   {\begin{itemize}}
\newcommand{\ei}   {\end  {itemize}}

\def\eq#1\en{\begin{equation}#1\end{equation}}  
\def\eqsplit#1\ensplit{
	\begin{equation}\begin{split}#1\end{split}\end{equation}
	}
\def\eqalign#1\enalign{
	\begin{align}#1\end{align}
	}
\def\eqmul#1\enmul{
	\begin{multline}#1\end{multline}
	}
\newcommand{\eqarrstar} {\begin{eqnarray*}} 
\newcommand{\enarrstar} {\end{eqnarray*}} 
\newcommand{\eqarray}   {\begin{eqnarray}} 
\newcommand{\enarray}   {\end{eqnarray}}

\newcommand{\lbeq}[1]  {\label{e:#1}}
\newcommand{\refeq}[1] {\eqref{e:#1}}    

\makeatletter
\newcommand{\labelcounter}[2]{{%
	\stepcounter{#1}
	\protected@write\@auxout{}%
	{\string\newlabel{#2}{{\csname the#1\endcsname}{\thepage}}}%
	{\ref{#2}}
	}}
\makeatother
\newcommand{\sss}   { \scriptscriptstyle } 

\newcommand{\Ncal}   {\mathcal{N}}

\newcommand{\Tcal}   {\mathcal{T}}

\newcommand{\Wcal}   {\mathcal{W}}

\newcommand{\spose}[1] {{\hbox to 0pt{#1\hss}} }
\newcommand{\ltapprox} {\mathrel{\spose{\lower 3pt\hbox{$\mathchar"218$}}
 \raise 2.0pt\hbox{$\mathchar"13C$}}}
\newcommand{\gtapprox} {\mathrel{\spose{\lower 3pt\hbox{$\mathchar"218$}}
 \raise 2.0pt\hbox{$\mathchar"13E$}}}

\newcommand{\dist} {{  \rm dist }}

\newcommand{\expec}[1]	{\left \langle #1 \right \rangle}

\usepackage[usenames]{color}

\newcommand{\ch}[1]{#1}

\newtheorem{theorem}{Theorem}[section]

\newtheorem{lemma}[theorem]{Lemma}
\newtheorem{corollary}[theorem]{Corollary}

\newtheorem{proposition}[theorem]{Proposition}
\newtheorem{conjecture}[theorem]{Conjecture}

\newtheorem{rem}[theorem]{Remark}
\newcommand{\eqa}{\begin{eqnarray}}
\newcommand{\ena}{\end{eqnarray}}
\newcommand{\nn}{\nonumber}

\newcommand{\whp}{{\bf whp}}

\newcommand{\prob}{\mathbb P}
\renewcommand{\expec}{\mathbb E}
\newcommand{\Var}{{\rm Var}}
\newcommand{\Cov}{{\rm Cov}}

\newcommand{\eqn}[1]{\begin{equation} #1 \end{equation}}
\newcommand{\eqan}[1]{\eqalign #1 \enalign}

\newcommand {\vep}{\varepsilon}

\newcommand {\sG}{{\sss G}}

\def\1{{\mathchoice {1\mskip-4mu\mathrm l}      
{1\mskip-4mu\mathrm l}
{1\mskip-4.5mu\mathrm l} {1\mskip-5mu\mathrm l}}}

\numberwithin{equation}{section}

\renewcommand{\to}      {\rightarrow}

\newcommand{\ben}{\begin{enumerate}}
\newcommand{\een}{\end{enumerate}}

\newcommand{\N}{{\mathbb N}}

\newcommand{\diam}{{\rm diam}}

\newcommand{\Core}{{\rm Core}_{t}}

\newcommand{\Inner}{{\rm Inner}_{t}}
\renewcommand{\Ncal}[1]{{\cal N}^{\sss(#1)}}
\renewcommand{\Wcal}[1]{u_{#1}}

\renewcommand{\Tcal}[2]{{\cal T}^{\sss(#1)}_{#2}}

\def\eqalign#1\enalign{
    \begin{align}#1\end{align}
    }

\title  {
        Diameters in preferential attachment models
        }

\author
{
Sander Dommers and Remco van der Hofstad
\footnote{
Eindhoven University of Technology, Department of Mathematics and Computer Science,
P.O.\ Box 513, 5600 MB Eindhoven, The Netherlands.
E-mail: {\tt S.Dommers@tue.nl, rhofstad@win.tue.nl}}\\
Gerard Hooghiemstra\footnote{Delft University of Technology,
Electrical Engineering, Mathematics and Computer Science, P.O. Box
5031, 2600 GA Delft, The Netherlands. E-mail: {\tt
G.Hooghiemstra@ewi.tudelft.nl} }}

\date{\today}

\begin{document}

\maketitle

\begin{abstract}
In this paper, we investigate the diameter in preferential
attachment (PA-) models, thus quantifying the statement that these
models are small worlds. The models studied here are such that
edges are attached to older vertices proportional to the degree
plus a constant, i.e., we consider {\it affine} PA-models. There
is a substantial amount of literature proving that, quite
generally, PA-graphs possess power-law degree sequences with a
power-law exponent $\tau>2$.

We prove that the diameter of the
PA-model is bounded above by a constant times $\log{t}$, where $t$
is the size of the graph. When the power-law exponent $\tau$
exceeds $3$, then we prove that $\log{t}$ is the right order, by proving a lower bound of this order,
both for the diameter as well as for the \ch{typical}
distance. This shows that, for $\tau>3$, distances are of the order
$\log{t}$. For $\tau\in (2,3)$, we improve the upper bound to a
constant times $\log\log{t}$, and prove a lower bound of the same
order for the diameter. Unfortunately, this proof does not extend
to \ch{typical} distances. These results do show that the diameter
is of order $\log\log{t}$.

These bounds partially prove predictions by physicists
that the \ch{typical} distance in PA-graphs are similar to the ones in
other scale-free random graphs, such as the configuration model and
various inhomogeneous random graph models, where \ch{typical} distances have been
shown to be of order $\log\log{t}$ when $\tau\in (2,3)$, and of
order $\log{t}$ when $\tau>3$.
\end{abstract}


\section{Introduction}
\label{sec-intro}
In the past decade, many examples have been found of real-world
complex networks that are {\it small worlds} and
{\it scale free}. The small-world phenomenon
states that distances in networks
are small. The scale-free phenomenon states
that the degree sequences in these networks
satisfy a power law. See \cite{AlbBar02, DorMen02, Newm03a}
for reviews on complex networks, and \cite{Bara02} for
a more expository account. Thus,
these complex networks are not at all like
classical random graphs (see \cite{AloSpe00, Boll01, JanLucRuc00}
and the references therein), particularly
since the classical models do not have power-law
degrees. As a result, these empirical findings
have ignited enormous research on random graph models that do obey power-law degree
sequences. See \cite{BolJanRio07} for the most general inhomogeneous random graph models,
as well as a review of the models under investigation.
Extensive discussions of various scale-free random graph models
are given in \cite{ChuLu06c, Durr06}.

While these models have power-law degree sequences,
they do not explain {\it why} many complex networks
are scale free. A possible explanation was given by
Barab\'asi and Albert \cite{BarAlb99} by a phenomenon called
{\it preferential attachment} (PA). Preferential attachment
models the growth of the network in such a way that
new vertices are more likely to add their edges
to already present vertices having a high degree.
For example, in a social network, a newcomer
is more likely to get to know a person who is
socially active, and, therefore, already has a high
number of acquaintances (high degree). Interestingly, PA-models with so-called \emph{affine} PA rules
have power-law degree sequences, and, therefore,
preferential attachment offers a convincing
explanation why many real-world networks possess this property.
There is a large amount of literature studying such models. See e.g.\
\cite{AieChuLu02, BolBorChaRio03, BolRio03b, BolRio03a,
BolRio04a, BolRio04b, BolRioSpeTus01, CooFri03}
and the references therein. The literature primarily
focusses on three main questions. The first key question for PA-models
is to prove that such random graphs are indeed {\it scale free}
\cite{AieChuLu02, BolBorChaRio03, BolRio03b, BolRio03a,
 BolRioSpeTus01, CooFri03}, by proving that their degree sequence
 indeed obeys a power law with a certain power-law exponent
 $\tau>2$. The second key question for PA-models
is their {\it vulnerability}, for example to deliberate
attack \cite{BolRio03a} or to the spread of a disease \cite{BerBolBorChaRio03}.
The third key question for PA-models is to show that the resulting
models are \emph{small worlds} by investigating the distances in them.
See in particular \cite{BolRio04b} for a result on the diameter
for a PA-model with power-law exponent $\tau=3$.
In non-rigorous work, it is often suggested that many
of the scale-free models, such as the configuration model,
the inhomogeneous random graph models in \cite{BolJanRio07} and the PA-models, have
similar properties for their distances. Distances
in the configuration model have been shown to depend on
the number of finite moments of the degree distribution. Similar results are
true for the so-called rank-1 inhomogeneous random graph (see e.g. \cite{ChuLu02a,ChuLu03,EskHofHoo06,NorRei06}).
The natural question is, therefore, whether the same applies to
preferential attachment models. This is the main goal of the present paper, in
which we investigate the diameter of scale-free PA-models.

The remainder of this section is organized as follows. We first introduce the models that
we will investigate in this paper. Then we give the main results and conclude with a discussion
of universality in power-law random graphs.

In this paper, we investigate the diameter in some PA-models. The models that we
investigate produce a {\it graph sequence} or {\it graph process}
$\{G_{m,\delta}(t)\}$, which, for fixed  $t\ge 1$ or $t \geq 2$, yields a graph
with $t$ vertices and $mt$ edges for some given integer
$m\geq 1$. In the sequel, we shall denote the vertices
of $G_{m,\delta}(t)$ by $1^{\sss(m)},\ldots,t^{\sss(m)}$. When $m$ is clear from the context,
we will leave out the superscript and write $[t]\equiv \{1,2,\ldots, t\}$.
We shall consider three slight variations of the
PA-model, which we shall denote by models (a), (b) and (c), respectively.

\begin{itemize}
\item[(a)]
The first model is an extension of the Barab\'asi-Albert model
formulated rigorously in  \cite{BolRioSpeTus01}.
We start with $G_{1,\delta}(1)$ consisting of a single vertex with
a single self-loop. We denote the degree of vertex $i^{\sss(1)}$ at time $t$ by $D_{i^{\sss (1)}}(t)$, where
a self-loop increases the degree by 2.\\
Then, for $m=1$, and conditionally on $G_{1,\delta}(t)$, the growth rule to obtain $G_{1,\delta}(t+1)$ is
as follows. We add a single vertex $(t+1)^{\sss(1)}$ having a single edge. This edge
is connected to a second end point, which is equal to $(t+1)^{\sss(1)}$ with probability proportional to $1+\delta$, and to a vertex $i^{\sss(1)}\in G_{1,\delta}(t)$
with probability proportional to $D_{i^{\sss(1)}}(t)+\delta$, where
$\delta\geq -1$ is a parameter of the model. Thus,
    \eqn{
    \label{growthrulePA}
    \prob\Big((t+1)^{\sss(1)}\rightarrow i^{\sss(1)}\big|G_{1,\delta}(t)\Big)=\left\{
    \begin{array}{lll}
    &\frac{1+\delta}{t(2+\delta)+(1+\delta)}, &\text{ for }i=t+1,\\
    &\frac{D_{i^{\sss(1)}}(t)+\delta}{t(2+\delta)+(1+\delta)}, &\text{ for }i\in[t].
    \end{array}
    \right.
    }
The model with integer $m>1$, is defined in terms of the model for $m=1$ as follows.
We start with $G_{1,\delta'}(mt),$ with $\delta'= \delta/m \ge -1$. Then we identify the vertices $1^{\sss(1)},2^{\sss(1)} \ldots, m^{\sss(1)}$
in $G_{1,\delta}(mt)$ to be vertex $1^{\sss(m)}$ in $G_{m,\delta}(t)$, and for $1< j \leq t$, the vertices
$((j-1)m+1)^{\sss(1)}, \ldots, (jm)^{\sss(1)}$ in $G_{1,\delta'}(mt)$ to be vertex $j^{\sss(m)}$ in $G_{m,\delta}(t)$; in particular
the degree $D_{j^{\sss(m)}}(t)$ of vertex $j^{\sss(m)}$ in $G_{m,\delta}(t)$ is equal to the sum of the degrees of the vertices
$((j-1)m+1)^{\sss(1)}, \ldots, (jm)^{\sss(1)}$ in $G_{1,\delta'}(mt)$.
This defines the model for integer $m\geq 1$. Observe that the range of $\delta$ is $[-m,\infty)$.\\
The resulting graph
$G_{m,\delta}(t)$ has precisely $mt$ edges and
$t$ vertices at time $t$, but is not
necessarily connected. For $\delta=0$ we obtain the original model
studied in \cite{BolRioSpeTus01}, and further studied in
\cite{BolRio03a, BolRio04a, BolRio04b}. The extension to $\delta\neq 0$
is crucial in our setting, as we shall explain in more detail below.

\item[(b)] The second model is identical to the one
above, apart from the fact that no self-loops are allowed for $m=1$.
We start again with the definition for $m=1$. To prevent a self-loop
in the first step, we let $G_{1,\delta}(1)$ undefined, and start
from $G_{1,\delta}(2)$, which is defined by the vertices $1^{\sss(1)}$ and $2^{\sss(1)}$ joined together
by 2 edges. Then, for $t \geq 2,$
we define, conditionally on $G_{1,\delta}(t)$, the growth rule to obtain $G_{1,\delta}(t+1)$
as follows. For $\delta \geq -1$,
\eq
\label{growthrulePA(b)}
    \prob\Big((t+1)^{\sss(1)}\rightarrow i^{\sss(1)}\big|G_{1,\delta}(t)\Big)=
    \frac{D_{i^{\sss(1)}}(t)+\delta}{t(2+\delta)}, \qquad\text{ for }i\in[t].
\en
The model with $m>1$ is again defined in terms of the model for $m=1$,
in precisely the same way as in model (a). This model is studied in detail in
\cite{Durr06}, and the model with $m=1$ corresponds to scale-free trees as
studied in e.g.\ \cite{BolRio04c, Mori02, Mori05, Pitt94}.

\item[(c)]
In the third model, and conditionally on $G_{m,\delta}(t)$, the end
points of each of the $m$ edges of vertex $t+1$, are chosen
{\it independently}, and are equal to a vertex $i^{\sss(m)}\in G_{m,\delta}(t),$
with probability proportionally to  $D_{i^{\sss(m)}}(t)+\delta$, where
$\delta\geq -m$.
We start again from $G_{m,\delta}(2)$, with the vertices $1^{\sss(m)}$ and $2^{\sss(m)}$ joined together by $2m\, , m\ge 1,$
edges. Since the end point of the edges are chosen independently we can give the definition
of $\{G_{m,\delta}(t)\}_{t\ge 2}$, for $m\ge 1$, in one step. For $1\le j\le m$,
    \eqn{
    \label{growthrulePAc}
    \prob\Big(j^{\rm th} \text{ edge of }(t+1)^{\sss(m)} \text{ is connected to }
    i^{\sss(m)}\big|G_{m,\delta}(t)\Big)=
    \frac{D_{i^{\sss(m)}}(t)+\delta}{t(2m+\delta)}, \qquad \text{ for }i\in[t].
    }
In this model, as is the case in model (b), the graph $G_{m,\delta}(t)$ is a connected
random graph with precisely $t$ vertices and $mt$ edges. This model was
studied in \cite{DeiEskHofHoo06, Jord06}.
\end{itemize}

\begin{rem} In models (a) and (b) for $m>1$, the choice of
$\delta'=\delta/m$ is such that in the
resulting graph $G_{m,\delta}(t)$, where $m$ vertices in $G_1(mt)$ are
grouped together to a single vertex in $G_{m,\delta}(t)$,
the end points of the added edges are chosen according
to the degree plus the constant $\delta$.
\end{rem}

\begin{rem} For $m=1$, the models (b) and (c) are the same.
This fact will be used later on.
\end{rem}

The growth rules in \eqref{growthrulePA}--\eqref{growthrulePAc}
are indeed such that vertices with high degree are more likely
to attract edges of new vertices. One would expect the models (a)--(c) to behave quite
similarly, as is known rigorously for the scale-free behaviour,
where the asymptotic degree distribution is known to be equal
in models (a)--(c).
As it turns out, the affine PA mechanism in
\eqref{growthrulePA}--\eqref{growthrulePAc} gives rise
to power-law degree sequences. Indeed, in \cite{DeiEskHofHoo06},
it was proved that for model (c), the degree sequence is close to a
power law with exponent $\tau=3+\delta/m$. For model (a)
and $\delta=0$, this was proved in \cite{BolRioSpeTus01}, while in
\cite{CooFri03}, power-law degree sequences for
PA-models with affine PA mechanisms are proved in rather
large generality. We see that, by varying the
parameters $m\geq 1, \delta>-m$, we can obtain
any power-law exponent $\tau>2$, which is the reason for
introducing the parameter $\delta$ in
\eqref{growthrulePA}--\eqref{growthrulePAc}. However, there is
no intrinsic reason for the {\it affine} PA mechanism.
For results on PA-models in the non-affine case,
see e.g., \cite{OliSpe05, RudTotVal07}. In general,
such models do {\it not} produce power laws.

The goal in this paper is
to study the diameter in the above
models, as a first step towards the study of distances
in PA-models and the verification of the prediction
that distances behave similarly in various scale-free
random models (see also Section \ref{sec-disc} below).
In the following section, we describe our precise results.

\subsection{Bounds on the diameter in preferential attachment models}
\label{sec-resdiam}
In this section, we present the diameter results for the PA-models (a)--(c).
The \emph{diameter} of a graph $G$ is defined as
    \begin{equation}
    \diam(G) = \max_{i,j\in G} \{\dist_G (i,j) | \dist_G (i,j) < \infty\},
    \end{equation}
where $\dist_G (i,j)$ denotes the graph distance between vertices $i,j\in G$.
We prove that, for all $\delta>-m,$ the diameter of
$G_{m,\delta}(t)$ is bounded by a constant times $\log{t}$.
When $\delta=0$, we adapt the argument
in \cite{BolRio04b} to prove that the diameter is bounded from below
by $(1-\vep)\frac{\log{t}}{\log\log{t}}$. For $\delta>0$, this lower bound is
improved to a constant times $\log{t}$, while, for $\delta<0$, we
prove that the diameter is bounded above and below by a constant times
$\log\log{t}$. This establishes a phase transition for the diameter
of PA-models when $\delta$ changes sign. We now state the precise results,
which shall all hold for each of the models (a)--(c) simultaneously.
In the results below, for a sequence of events $\{A_t\}_{t\ge 1}$, we write
that $A_t$ occurs {\it with high probability} (\whp) when $\lim_{t\to \infty}
\prob(A_t)=1$.

\begin{theorem}[A $\log{t}$ upper bound on the diameter]
\label{thm-diamPA}
Fix $m\geq 1$ and $\delta>-m$. Then,
there exists a constant $c_1=c_1(m,\delta)>0$ such that
\whp, the diameter of $G_{m,\delta}(t)$ is at most $c_1\log t$.
\end{theorem}

When $m=1$, so that the graphs are in fact trees,
there is a sharper result proved by Pittel \cite{Pitt94},
which, in particular, implies Theorem \ref{thm-diamPA}
for model (b). In this case, Pittel shows that the height of the tree, which
is equal to the maximal graph distance between vertex 1 and
any of the other vertices, grows like $\frac{1+\delta}
{\gamma (2+\delta)}\log{t}(1+o(1))$,
where $\gamma$ solves the equation
    \eqn{
    \lbeq{gammaeq}
    \gamma +(1+\delta)(1+\log{\gamma})=0.
    }
This proves that the diameter is at least as large, and
suggests that the diameter has size $\frac{2(1+\delta)}
{\gamma (2+\delta)}\log{t}(1+o(1))$. Scale-free trees have received
substantial attention in the literature, we refer
to \cite{BolRio04c, Pitt94} and the references therein.
It is not hard to see that a similar result as proved in
\cite{Pitt94} also follows for models (a) and (c). This is proved
when $\delta=0$ in \cite{BolRio04c}, where it is shown
that the diameter in model (a) has size $\gamma^{-1} \log{t}$,
where $\gamma$ is the
solution of \refeq{gammaeq} when $\delta=0$. Thus, we see that
the $\log{t}$ upper bound in Theorem \ref{thm-diamPA} is sharp,
at least for $m=1$.

It is not hard to extend the
upper bound to $m\geq 2$. In particular, for model (b), the upper bound
for $m\geq 2$ immediately follows from the upper bound for $m=1$.
For models (a) and (c), the extension is not as trivial, but
the proof is fairly straightforward, and will be omitted here.
To see an implication of \cite{Pitt94} for model (a),
we note that $C_t$, the number of connected
components of $G_{1,\delta}(t)$ in model (a), has distribution $C_t=1+I_2+\cdots+I_t$,
where $I_i$ is the indicator that the $i^{\rm th}$ edge connects to itself,
so that $\{I_i\}_{i=2}^t$ are independent indicator variables with
    \eqn{
    \prob(I_i=1)=\frac{1+\delta}{(2+\delta)(i-1)+1+\delta}.
    }
As a result, $C_t/\log{t}$ converges in probability to $(1+\delta)/(2+\delta)<1$,
so that {\bf whp}  there exists a largest connected component of size at least
$t/\log{t}$. \ch{Conditionally on having size $s_t$, t}he law of any connected component in model (a) is
{\it equal} in distribution to the law of the graph $G_{1,\delta}(s_t+1)$ in model (b), apart
from the fact that the vertices 1 and 2 in $G_{1,\delta}(s_t+1)$ are identified
(thus creating a \ch{double} self-loop) \ch{and the vertices are relabeled by order of appearance. In particular, conditionally on having size $s_t$, the law of the {\em diameter} of the connected component in model (a) equals that of $G_{1,\delta}(s_t+1)$ in model (b).} This close connection between the
two models allows one to transfer results for model (b) to model (a)
when $m=1$.

\begin{theorem}[A $\log{t}$ lower bound on the diameter for $\delta>0$]
\label{thm-diamPAlb}
Fix $m\geq 1$ and $\delta>0$. Then,
there exists $c_2=c_2(m,\delta)>0$, such that
\whp, the diameter of $G_{m,\delta}(t)$ is at least $c_2\log{t}$.
\end{theorem}

Theorems \ref{thm-diamPA}--\ref{thm-diamPAlb} imply that,
for $\delta>0$ and {\bf whp}, $\diam(G_{m,\delta}(t))=\Theta(\log{t})$.
Theorems \ref{thm-diamPA}--\ref{thm-diamPAlb} indicate that distances in
PA-models are similar to the ones in other scale-free models for $\tau>3$.
We shall discuss this analogy in more detail below.
As we shall see in Section \ref{sec-PA>0},
the proof of Theorem \ref{thm-diamPAlb} also reveals that,
{\bf whp}, the \ch{typical} distance in $G_{m,\delta}(t)$,
which is the distance between two uniformly chosen connected vertices
in the graph, is bounded from below by $c_2 \log{t}$.

We conjecture that, for $\delta>0$, a limit result holds for the
constant in front of the $\log{t}$. In its statement, we write $\dist_{G}(v_1,v_2)$
for the graph distance in the graph $G$ between two vertices $v_1,v_2\in [t]$.
Then, the \ch{typical} distance in a graph $G$ is defined by $\dist_{G}(V_1,V_2)$
where $V_1,V_2\in [t]$ are two uniformly chosen independent vertices.

\begin{conjecture}[Convergence in probability for $\delta>0$]
\label{conj-diamPAlb}
Fix $m\geq 1$ and $\delta>0$. Then,
the diameter
$\diam(G_{m,\delta}(t))/\log{t}$ and the \ch{typical} distance
$\dist_{G}(V_1,V_2)/\log{t}$ converge in probability
to positive and different constants.
\end{conjecture}

We now turn to the case where $\delta \in (-m,0)$ and hence $\tau=3+\delta/m \in (2,3)$:

\begin{theorem}[A $\log\log{t}$ upper bound on the diameter for $\delta<0$]
\label{thm-diametersec} Fix $m\geq 2$ and assume that $\delta \in
(-m,0)$. Then, for every $\sigma>1/(3-\tau)$ and with
    \eq
    \lbeq{CG-def}
    C_{\sG}=\frac{4}{|\log{(\tau-2)}|}+\frac{4\sigma}{\log{m}},
    \en
the diameter of $G_{m,\delta}(t)$ is, {\bf whp}, bounded above by $C_{\sG} \log\log{t}$,
as $t\rightarrow \infty$.
\end{theorem}

In this result, we do not obtain a sharp result in
terms of the constant. However, the proof
suggests that for most pairs of vertices the distance should be equal
to $\frac{4}{|\log{(\tau-2)}|}\log{\log {t}}(1+o(1))$.
When $m=1$, Theorem \ref{thm-diametersec} does {\it not} hold
(see the discussion below Theorem \ref{thm-diamPA}).

We next discuss the lower bound on the diameter for $\delta\in(-m,0)$:

\begin{theorem}[A $\log\log{t}$ lower bound on the diameter]
\label{thmlowerb}
Fix $m\geq2$ and $\delta>-m$. Then, the diameter of $G_{m,\delta}(t)$ is, \whp,
bounded below by $\frac{\varepsilon}{\log{m}}\log{\log{t}}$, for all $\varepsilon \in (0,1)$.
\end{theorem}

Unfortunately, the proof of Theorem \ref{thmlowerb} does not allow for an extension
to \ch{typical} distances, and, thus, we have no matching lower bound for this.
We finally conjecture that, for $\delta \in(-m,0)$, a limit results holds for the
constant in front of the $\log\log{t}$:

\begin{conjecture}[Convergence in probability for $\delta<0$]
\label{conj-diamPAlb}
Fix $m\geq 2$ and $\delta\in(m,0)$. Then, the diameter
$\diam(G_{m,\delta}(t))/\log{\log{t}}$ and the \ch{typical} distance
$\dist_{G}(V_1,V_2)/\log{\log{t}}$ converge in probability
to positive and different constants.
\end{conjecture}

\subsection{Discussion of universality of distances in power-law random graphs}
\label{sec-disc}
Theorems \ref{thm-diamPA}--\ref{thmlowerb} prove that the diameter in
PA-models with a power-law degree sequence denoted by $\tau$ undergoes a
phase transition as $\tau$ changes from $\tau\in (2,3)$ to $\tau>3$.
The results identify the order of growth of the diameter of three related
models of affine PA models as the size of the graph $t$ tends to infinity.
We do not obtain the right constants. For the \ch{typical} distance,
we obtain a similar phase transition, and again the results identify
the correct asymptotics for $\tau>3$, but, for $\tau\in (2,3)$
we miss a matching lower bound.

In non-rigorous work, it is often suggested that the distances
are similarly behaved in the various scale-free random graph models,
such as the configuration model or various models with conditional
independence of edges as in \cite{BolJanRio07}. For power-law
random graphs, this informal statement can be made precise
by conjecturing that distances have the same leading order
growth in graphs with the same power-law degree exponent.
This, however, is not correct for the {\it diameter}
of such power-law random graphs, since the diameter depends sensitively
on the details of the graph, such as the proportion of vertices with
degrees 1 and 2. See \cite{FerRam04} and \cite{HofHooZna07} for
results showing that for the configuration model with power-law degree
exponent $\tau\in (2,3)$, the diameter
can be of order $\log{t}$ or of order $\log{\log{t}}$ depending
on the proportion of vertices with degrees 1 and 2,
where $t$ is the size of the graph.
Similarly, in inhomogeneous random graphs with power-law
degree exponent $\tau\in (2,3)$ the diameter is always of order
$\log{t}$ (see e.g.\ \cite{BolJanRio07}), while the
\ch{typical} distances can be of order
$\log\log{t}$ (see e.g.\ \cite{ChuLu02a, ChuLu03}).
Thus, we shall interpret the physicists' prediction
by conjecturing that the leading order growth of the
{\it \ch{typical}} distances of various power-law random graphs
depends only on the power-law degree exponent $\tau\in (2,3)$.

The results on distances are most complete for the configuration model (CM),
see e.g.\ \cite{EskHofHooZna04, FerRam04, HHV05, HofHooZna04a, ReiNor04}.
In the CM, there
are various cases depending on the tails of the degree distribution.
When the degrees have infinite mean, then \ch{typical}
distances are bounded \cite{EskHofHooZna04},
when the degrees have finite mean but infinite variance, \ch{typical} distances grow
proportionally to $\log\log{t}$ \cite{HofHooZna04a, ReiNor04}, where $t$ is the size
of the graph, while, for finite
variance degrees, the \ch{typical} distances grow proportionally to $\log{t}$ \cite{HHV05}.
Similar results for models with conditionally independent edges exist,
see e.g. \cite{BolJanRio07, ChuLu02a,EskHofHoo06, NorRei06}, but
particularly in the regime $\tau\in (2,3)$, the results are not that
strong.
Thus, for these classes of models, distances are quite well understood.
If the distances in PA-models are similar
to the ones in e.g.\ the CM, then we should have
that the distances are of order $\log{t}$ when $\tau>3$, i.e., $\delta>0$,
while they should be of order $\log{\log{t}}$ when $\tau\in (2,3)$,
i.e., for $\delta<0$. In PA-models with
a linear growth of the number of edges, infinite mean degrees
cannot arise, which explains why $\tau>2$ for PA-models.
An attempt in the direction of creating PA-models with power-law
exponent $\tau\in (1,2)$ can be found in \cite{DeiEskHofHoo06},
where a preferential attachment model is presented in which a
{\it random} number of edges per new vertex is added.
In this model, it is shown that the degrees again obey
a power law with exponent equal to $\tau=\min\{3+\frac{\delta}{\mu}, \tau_w\}$,
where $\tau_w$ is the power-law exponent for the number of edges
added and $\mu\leq \infty$ the expected number of added edges per vertex.
Thus, when $\tau_w\in (1,2)$, infinite mean degrees can arise.
This model is further studied in \cite{Bham07}, where a wealth
of results for various PA-models can be found.

There are few results on distances in PA-models.
In \cite{BolRio04b}, it was proved that in model (a)
and for $\delta=0$, for which $\tau=3$, the diameter of the graph
of size $t$ is
equal to $\frac{\log{t}}{\log\log{t}}(1+o(1))$. Unfortunately, the matching result
for the CM has not been proved, so that this does not
allow us to verify whether the models have similar distances.
The results stated above substantiate the physicists' prediction,
since, for $\delta>0$ for which $\tau\in (3,\infty)$, the \ch{typical}
distances are of order $\log{t}$, while, for $\delta<0$, for which
$\tau\in (2,3)$, they are bounded above by $\log\log{t}$.
A related result on PA-models in the spirit
of \cite{CooFri03} can be found in \cite{ChuLu04b}, where a
similar phase transition as in this paper is proved,
in the case where the number of edges grows at least
$(\log{t})^{1+\varepsilon}$ times as fast as the number
of vertices.

It would be of interest to improve the bounds presented in this
paper up to the constant in front of the $\log{t}$ and
$\log\log{t}$, respectively. Due to the dynamical nature of
PA-models, this is more involved for PA-models than it is for
static models such as the CM and inhomogeneous
random graphs.
\medskip

This paper is organized as follows. In Section \ref{sec-logLB},
we prove the $\log{t}$ lower bound for the diameter
stated in Theorem \ref{thm-diamPAlb}. In
Section \ref{sec-loglog} and Section \ref{lowerbound}, we prove the $\log\log{t}$ upper
bound and the $\log\log{t}$ lower bound, on the diameter for $\delta<0$, of Theorem \ref{thm-diametersec}
and Theorem \ref{thmlowerb}, respectively.


\section{A $\log$ lower bound on the diameter for $\delta> 0$: Proof of Theorem \ref{thm-diamPAlb}}
\label{sec-logLB}
In this section, we prove Theorem \ref{thm-diamPAlb} by extending the
argument in \cite{BolRio04b} from $\delta=0$ to $\delta>0$.
We shall also extend the lower bound for $\delta=0$ to models
(b) and (c).

For model (c), denote by
    \eq
    \label{gerard7}
    \{g(t,j)=s\}, \quad 1\le j\le m,
    \en
the event that at time $t$ the $j^{\rm th}$ edge of vertex $t$
is attached to the earlier vertex $s<t$. For models (a) and (b), this event means that
in $\{G_{1,\delta'}(mt)\}$ the edge from vertex $m(t-1)+j$ is attached to one of the vertices $m(s-1)+1,\ldots,ms$.
It is a direct consequence of the definition of PA-models that the event
\eqref{gerard7} increases the preference for vertex $s$, and hence decreases (in a relative way)
the preference for the vertices $u, \, 1\le u \le t,\, u\neq s$. It should be intuitively
clear that another way of expressing this effect is to say that, for different
$s_1\neq s_2$, the events $\{g(t_1,j_1)=s_1\}$ and $\{g(t_2,j_2)=s_2\}$  are negatively correlated.
In order to state such a result, we introduce some notation.
For integer $n_s\ge 1$ and $i=1, \ldots, n_s$, we denote by
    \eqn{
    \lbeq{Es-def}
    E_s=\bigcap_{i=1}^{n_s} \big\{g(t_i,j_i)=s\big\},
    }
the event that at time $t_i$ the
$j_i^{\rm th}$ edge of vertex $t_i$ is attached to the earlier vertex
$s$\ch{, for all $i=1,\ldots,n_s$}. We will start by proving that for each $k\geq 1$
and all possible choices of $t_i,j_i$, the events $E_s$, for different $s$, are negatively correlated:

\begin{lemma}[Negative correlation of attachment events] For distinct $s_1,s_2,\ldots,s_k$,
\label{lem-negcorE}
\eq
    \lbeq{negcorE}
    \prob\Big(\bigcap_{i=1}^k E_{s_i}\Big)\leq \prod_{i=1}^k \prob(E_{s_i}).
    \en
\end{lemma}

\proof We will use induction on the largest edge number present
in the events $E_s$. Here, for an event $\{g(t,j)=s\}$,
we let the edge number be $m(t-1)+j$, which is the order
of the edge when we consider the edges as being attached
in sequence. The induction hypothesis is that
\refeq{negcorE} holds for all $k$ and all choices of
$t_i,j_i$ such that $\max_{i,s}
m(t_i-1)+j_i \leq e$, where induction is
performed with respect to $e$. \ch{We initialize the induction for $e=m$ in models (a) and (b) and for $e=2m$ in model (c). We note that
for this choice of $e$, the induction hypothesis holds trivially, since everything is deterministic.}
This initializes the induction.

To advance the induction, we assume that \refeq{negcorE} holds for all $k$ and all choices of
$t_i,j_i$ such that $\max_{i,s}
m(t_i-1)+j_i \leq e-1$.
Clearly, for $k$ and $t_i,j_i$ such that
$\max_{i,s}
m(t_i-1)+j_i \leq e-1$, the
bound follows from the induction hypothesis, so we may
restrict attention to the case that
$\max_{i,s}
m(t_i-1)+j_i = e$.
We note that there is a unique choice of
$t,j$ such that $m(t-1)+j=e$. In this case, there are again
two possibilities. Either there is exactly one choice of
$s$ and $t_i,j_i$ such that
$t_i=t, j_i=j$, or there are at least
{\it two} of such choices. In the latter case, we
immediately have that \ch{$\bigcap_{i=1}^k E_{s_i}=\varnothing,$}
since the $e^{\rm th}$ edge can only be connected
to a {\it unique} vertex. Hence, there is nothing to prove.
Thus, we are left to investigate the case where there
exists unique $s$ and $t_i,j_i$ such that
$t_i=t, j_i=j$. Denote by
    \eqn{
    E_s'=
    \bigcap_{i=1: (t_i, j_i)\neq (t,j)}^{n_s}
    \big\{g(t_i,j_i)=s\big\},
    }
the restriction of $E_s$ to the {\it other} edges.
Then we can write
    \eqn{
    \bigcap_{i=1}^k E_{s_i}
    =\big\{g(t,j)=s\big\}
    \cap E_s' \cap  \bigcap_{i=1: s_i\neq s}^k E_{s_i}.
    }
By construction, all the edge numbers of the events
in $E_s' \cap \bigcap_{i=1: s_i\neq s}^k E_{s_i}$ are at most $e-1$.
Thus, we obtain
    \eqn{
    \prob\Big(\bigcap_{i=1}^k E_{s_i}\Big)
    \leq \expec\Big[I[E_s' \cap  \bigcap_{i=1: s_i\neq s}^k E_{s_i}]
    \prob_{e-1}(g(t,j)=s)\Big],
    }
where $\prob_{e-1}$ denotes the conditional probability given the
edge attachments up to the $(e-1)^{\rm st}$ edge connection,
and where, for an event $A$, $I[A]$ denotes the indicator of $A$.

We now first treat model (c), for which we have that
    \eqn{
    \prob_{e-1}(g(t,j)=s)= \frac{D_s(t-1)+\delta}{(2m+\delta)(t-1)}.
    }
We wish to use the induction hypothesis. For this, we note that
    \eqn{
    \lbeq{degree-form}
    D_s(t-1)=m+\sum_{(t',j'): t'\leq t-1} I[g(t',j')=s].
    }
We note  that
each of the terms in \refeq{degree-form} has edge number
strictly smaller than $e$ and occurs with a non-negative
multiplicative constant. As a result, we may use
the induction hypothesis for each of these terms. Thus,
we obtain, using also $m+\delta\geq 0$, that,
\eqa
(2m+\delta)(t-1)\prob\Big(\bigcap_{i=1}^k E_{s_i}\Big)
&&\leq (m+\delta)    \prob(E_s')\prod_{i=1: s_i\neq s}^k \prob(E_{s_i})\nn\\
&&\qquad+    \sum_{(t',j'): t'\leq t-1}\prob(E_s'\cap \{g(t',j')=s\})
    \prod_{i=1: s_i\neq s}^k \prob(E_{s_i}).
\ena
We can recombine to obtain
    \eqn{
    \prob\Big(\bigcap_{i=1}^k E_{s_i}\Big)
    \leq \expec\Big[I[E_s']
    \frac{D_s(t-1)+\delta}{(2m+\delta)(t-1)}\Big]\prod_{i=1: s_i\neq s}^k \prob(E_{s_i}) ,
    }
and the advancement is completed when we note that
    \eqn{
    \expec\Big[I[E_s']
    \frac{D_s(t-1)+\delta}{(2m+\delta)(t-1)}\Big]
    =\prob(E_s).
    }
The proofs for models (a) and (b) are somewhat simpler, since
the events $E_{s_i}$ can be reformulated in terms of the graph
process $\{G_{1,\delta'}(t)\}_{t\geq 1}$.
\qed

We next give  the probabilities of $E_s$ when $n_s\le 2$; we omit
the proof, since it is a simple adaptation to that in \cite{BolRio04b}.
\begin{lemma}[Connections in PA-models]
\label{lemmaE}
There exist absolute constants $M_1,M_2$, such that (i) for each $1\le j\le m$, and $t>s$,
    \eq
    \label{enkelprob}
    \prob\Big( g(t,j)=s\Big)\leq \frac{M_1}{t^{1-a}s^a},
    \en
and (ii) for $t_2>t_1>s$, and any $1\le j_1,j_2\le m$,
    \eq
    \label{simkans}
    \prob\Big(g(t_1,j_1)=s,g(t_2,j_2)=s\Big)\leq \frac{M_2}{(t_1t_2)^{1-a}s^{2a}},
    \en
    where $a=\frac{m}{2m+\delta}$.
\end{lemma}


\vskip 0.5truecm
\noindent
We combine the results of Lemmas \ref{lem-negcorE} and
\ref{lemmaE} into the following corollary, yielding an upper bound for the
probability of the existence of a path. In its statement, we call a path
$\Gamma=(s_0,s_1, \ldots, s_l)$ {\em self-avoiding} when $s_i\neq s_j$ for all $0\leq i<j\leq l$.
We use the notation $x \wedge y = \min(x,y)$ and $x \vee y = \max(x,y)$.
Again, we omit the  proof (for details, see \cite{BolRio04b}).

\begin{corollary}[Path probabilities in PA-models]
Let $\Gamma=(s_0,s_{1},\ldots,s_l)$ be a
self-avoiding path of length $l$ consisting of the
$l+1$ unordered vertices $s_0,s_{1},\ldots,s_l$, then
there exists an absolute constant $C>0$ such that
\eq
\label{unordupb}
\prob\Big(\Gamma\in G_{m,\delta}(t)\Big)\leq (m^2C)^{l}\prod_{i=0}^{l-1}
\frac{1}{(s_i\wedge s_{i+1})^a {(s_i\vee s_{i+1})}^{1-a}}.
\en
\end{corollary}

\subsection{Lower bound on the diameter for $\delta=0$}
\label{sec-PA0}
It follows from (\ref{unordupb}) that for  $\delta=0$,
    \eq
    \label{pathupb2}
    \prob\Big(\Gamma\in G_{m,\delta}(t)\Big)
    \leq (m^2C)^{l}\prod_{i=0}^{l-1}
    \frac{1}{\sqrt{s_i s_{i+1}}} .
    \en
The further proof that \eqref{pathupb2} implies that for $\delta \ge 0$,
    \eqa
    \label{lowerbndL}
    L=\frac{\log t}{\log (3Cm^2 \log t)},
    \ena
is a lower bound for the diameter of $G_{m,\delta}(t)$, is identical to the proof of
\cite[Theorem 5, p.~14]{BolRio04b}, with $n$ replaced by $t$. This extends
the lower bound for $\delta=0$ for model (a) in \cite{BolRio04b} to models
(b)--(c).
\qed

\subsection{The lower bound on distances for $\delta>0$}
\label{sec-PA>0}
We next improve the bound in the previous section in the case
when $\delta>0$, in which case $a=m/(2m+\delta)<1/2$.
From the above discussion, we conclude that
    \eqn{
    \label{diamGmk}
    \prob\Big(\dist_{\sss G_{m,\delta}(t)}(1,t)=k\Big)
    \leq c^k \sum_{\vec s} \prod_{j=0}^{k-1}
    \frac{1}{(s_j\wedge s_{j+1})^a {(s_j\vee s_{j+1})}^{1-a}},
    }
where $c=m^2C$, and where the sum is over $\vec s=(s_0, \ldots, s_k)$ with
$s_k=t, s_0=1$, $s_l\geq 1$ for all $l=1, \ldots, k-1$ and $s_l\neq s_n$ for all $l\neq n$\ch{, since we may assume that our path $(s_0,\ldots,s_k)$ is self-avoiding}. Define
    \eqn{
    f_k(i,t)=\sum_{\vec s} \prod_{j=0}^{k-1}
    \frac{1}{(s_j\wedge s_{j+1})^a {(s_j\vee s_{j+1})}^{1-a}},
    }
where now the sum is over $\vec s=(s_0, \ldots, s_k)$ with
$s_k=t, s_0=i$, $s_l\geq 1$ for all $l=1, \ldots, k-1$ and
$s_l\neq s_n$ for all $l\neq n$, so that
    \eqn{
    \label{diamGmkfk}
    \prob\Big(\dist_{\sss G_{m,\delta}(t)}(i,t)=k\Big)
    \leq c^k f_k(i,t).
    }
We study the function $f_k(i,t)$ in the following lemma:

\begin{lemma}[A bound on $f_k$]
\label{lem-fk}
Fix $a<1/2$. Then, for every $b>a$ such that $a+b<1$, there exists
a $C_{a,b}>0$ such that, for every $1\leq i < t$ and all $k \geq 1$,
    \eq
    \label{fkbd}
    f_k(i,t)\leq \frac{C_{a,b}^k}{i^b t^{1-b}}.
    \en
\end{lemma}

\proof We prove the lemma using induction on $k\geq 1$. To initialize the
induction hypothesis, we note that, for $1\leq i<t$
and every $b\geq a$,
    \eq
    f_1(i,t)= \frac{1}{(i\wedge t)^a {(i\vee t)}^{1-a}}
    =\frac{1}{i^a t^{1-a}}
    =\frac{1}{t} \Big(\frac{t}{i}\Big)^a
    \leq \frac{1}{t} \Big(\frac{t}{i}\Big)^b
    =\frac{1}{i^bt^{1-b}}.
    \en
This initializes the induction hypothesis
as long as $C_{a,b}\geq 1$.

To advance the induction hypothesis, note that we have
the recursion relation
    \eq
    f_k(i,t) \ch{\leq} \sum_{s=1}^{i-1}
    \frac{1}{s^a {i}^{1-a}}
    f_{k-1}(s,t)
    +\sum_{s=i+1}^{\infty}
    \frac{1}{i^a {s}^{1-a}}
    f_{k-1}(s,t).
    \en
We now bound each of these two contributions, making use of
the induction hypothesis. For the first sum, we bound
    \eq
    \sum_{s=1}^{i-1}
    \frac{1}{s^a {i}^{1-a}}
    f_{k-1}(s,t)
    \leq C^{k-1}_{a,b} \sum_{s=1}^{i-1}\frac{1}{s^a {i}^{1-a}}\frac{1}{s^b t^{1-b}}
    =\frac{C^{k-1}_{a,b}}{i^{1-a} t^{1-b}} \sum_{s=1}^{i-1}\frac{1}{s^{a+b}}
    \leq \frac{1}{1-a-b} \frac{C^{k-1}_{a,b}}{i^{b} t^{1-b}},
    \en
since $a+b<1$. For the second sum, we bound
\ch{\begin{align}
    \sum_{s=i+1}^{\infty}
    \frac{1}{i^a {s}^{1-a}}
    f_{k-1}(s,t)
    &\leq C^{k-1}_{a,b} \sum_{s=i+1}^{t-1}\frac{1}{i^a {s}^{1-a}}\frac{1}{s^b t^{1-b}}
    + C^{k-1}_{a,b} \sum_{s=t+1}^{\infty}\frac{1}{i^a {s}^{1-a}}\frac{1}{t^b s^{1-b}} \\
    &= \frac{C^{k-1}_{a,b}}{i^{a}t^{1-b}} \sum_{s=i+1}^{t-1}\frac{1}{s^{1-a+b}}
    + \frac{C^{k-1}_{a,b}}{i^{a}t^b} \sum_{s=t+1}^{\infty}\frac{1}{{s}^{2-a-b}}
    \leq \frac{1}{b-a} \frac{C^{k-1}_{a,b}}{i^{b} t^{1-b}} + \frac{1}{1-a-b}\frac{C^{k-1}_{a,b}}{i^{b} t^{1-b}},\nonumber
    \end{align}
}
since $1+b-a>1$\ch{, $2-a-b>1, b>a$ and $(t/i)^a \leq (t/i)^b$}. We conclude that
    \eq
    f_k(i,t)\leq \frac{C^{k-1}_{a,b}}{i^{b} t^{1-b}} \Big(\frac{1}{b-a}+\frac{\ch{2}}{1-a-b}\Big)
    \leq \frac{C^{k}_{a,b}}{i^{b} t^{1-b}},
    \en
when
    \eq
    C_{a,b}=\frac{1}{b-a}+\frac{\ch{2}}{1-a-b}\geq 1.
    \en
This advances the induction hypothesis, and completes the proof.
\qed
\vskip0.5cm

\noindent
Using Lemma \ref{lem-fk} and \eqref{diamGmkfk}, we obtain that
    \eqn{
    \prob\Big(\dist_{\sss G_{m,\delta}(t)}(1,t)=k\Big)
    \leq \frac{(cC_{a,b})^k}{t^{1-b}}.
    }
As a result,
    \eqn{
    \prob\Big(\diam(G_{m,\delta}(t))\leq k\Big)
    \leq \prob\Big(\dist_{\sss G_{m,\delta}(t)}(1,t)\leq k\Big)
    \leq \frac{(cC_{a,b})^{k+1}}{t^{1-b}(cC_{a,b}-1)}=o(1),
    }
whenever $k\leq \frac{1-b}{\log{(cC_{a,b})}} \log{t}$. We conclude
that there exists $c_2=c_2(m,\delta)$ such that, with high probability
$\diam\Big(G_{m,\delta}(t)\Big)\geq c_2\log{t}$.
\qed
\vskip0.5cm

\noindent
We next extend the above discussion to {\it \ch{typical}} distances.

\begin{lemma}[\ch{Typical} distances for $\delta>0$]
Fix $m\geq 1$ and $\delta>0$. Let $H_t={\rm dist}_t(A_1,A_2)$
be the distance between two uniformly chosen vertices.
Then, for $c_2=c_2(m,\delta)>0$ sufficiently small,  {\bf whp},
$H_t\geq c_2\log t$.
\end{lemma}

\proof
For $c_2=c_2(m,\delta)>0$, define
    \eqn{
    \lbeq{aver-dist}
    B_t\equiv \#\big\{i,j\in [t]: i<j: \dist_{\sss G_{m,\delta}(t)}(i,j)\leq c_2\log t\big\}
    ,
    }
where $\#\{A\}$ denotes the cardinality of the set $A$.

By Lemma \ref{lem-fk}, with $K=\log{(cC_{a,b}\vee 2)}$
and $a<b<1-a$, and for all $1\leq i<j\leq t$,
    \eqn{
    \prob\Big(\dist_{\sss G_{m,\delta}(t)}(i,j)=k\Big)\leq c^k f_k(i,j)
    \leq
    \frac{e^{Kk}}{\ch{i^{b}j^{1-b}}}.
    }
As a result,
    \eqn{
    \prob\Big(\dist_{\sss G_{m,\delta}(t)}(i,j)\leq c_2\log t\Big)
    \leq \frac{t^{Kc_2}}{\ch{i^{b}j^{1-b}}}\frac{e^K}{e^K-1},
    }
and thus, using also \ch{$\sum_{i=1}^{j-1} i^{-b}\leq j^{1-b}/(1-b)$},
    \eqn{
    \lbeq{expecBt}
    \expec[B_t]\leq O(1)\sum_{1\leq i<j\leq t} \frac{t^{Kc_2}}{\ch{i^{b}j^{1-b}}}
    =O(t^{Kc_2+1}).
    }
It now suffices to note that
    \eqn{
    \prob(H_t\leq c_2\log t)
    =\expec\big[I[\dist_{\sss G_{m,\delta}(t)}(A_1,A_2)\leq c_2 \log t]\big]
    \ch{=} \frac{2\expec[B_t]+t}{t^2}=o(1),
    }
by \refeq{expecBt}, for every $c_2>0$ such that $Kc_2+1<2$.
\qed

Note that \eqref{lowerbndL} is also a lower bound on \ch{typical} distances in case $\delta=0$, which can be proved as above.

\section{A $\log\log$ upper bound on the diameter: Proof of
Theorem~\ref{thm-diametersec}}
\label{sec-loglog}
The proof of Theorem \ref{thm-diametersec} is divided into two
key steps. In the first, in Theorem \ref{prop-core},
we bound the diameter
of the {\it core} which consists of
the vertices with degree at least a certain power of $\log{t}$.
This argument is close in spirit to the argument in \cite{ChuLu02a} or \cite{ReiNor04}
used to prove bounds on the \ch{typical} distance for the
inhomogeneous random graph and the configuration
model, respectively, but substantial adaptations are necessary to deal with
preferential attachment.
After this, in Theorem \ref{prop-periphery}, we derive a bound
on the distance between vertices with a small degree and the core.
We start by defining and investigating the core of the PA-model.
In the sequel, it will be convenient to prove Theorem
\ref{thm-diametersec} for $2t$ rather than for $t$. Clearly, this
does not make any difference for the results. We make use of some technical results, stated in the appendix.

\subsection{The diameter of the core}
We recall that $\tau=3+\delta/m$, so that $-m<\delta<0$ corresponds to $\tau\in (2,3)$.
We take $\sigma>1/(3-\tau)=-m/\delta>1$ and define the {\it core}
$\Core$ to be
    \eq \label{defcore}
    \Core=\big\{i\in [t]: D_i(t)\geq (\log{t})^\sigma\big\},
    \en
i.e., all the vertices which at time $t$ have degree at least $(\log{t})^{\sigma}$.

For $A\subseteq [t]$, we write
    \eq
    \diam_t(A) = \max_{i,j\in A} \dist_{\sss G_{m,\delta}(t)}(i,j).
    \en
Then, $\diam_{2t}(\Core)$ is bounded in the following theorem:

\begin{theorem}[The diameter of the core]
\label{prop-core}
Fix $m\geq 2$ and $\delta\in (-m,0)$. For every $\sigma>1/(3-\tau)$, \whp,
    \eq
    \diam_{2t}(\Core)\leq (1+o(1)) \frac{4 \log\log{t}}{|\log{(\tau-2)}|}.
    \en
\end{theorem}

The proof of Theorem \ref{prop-core} is divided into several
smaller steps. We start by proving that the diameter $\diam_{2t}(\Inner)$, where
    \eq
    \label{Inner-def}
    \Inner =\big\{i\in [t]: D_i(t)\geq u_1\big\}, \quad \mbox{and where}
    \quad u_1= t^{\frac{1}{2(\tau-1)}} (\log{t})^{-\frac{1}{2}},
    \en
is, \whp, bounded. The choice of $u_1$ is a technical one:
$u_1$ is the largest value $l$ so that, \whp, the total degree
of vertices with degree exceeding $l$ can be bounded from below
by $tl^{2-\tau}$, see Lemma \ref{totdegree}. In Proposition \ref{prop-inner},
we will show that the diameter of $\Inner$ is bounded.
After this, we will show that
the distance from any vertex in the core $\Core$
to the inner core $\Inner$ can be bounded by
a fixed constant times $\log\log{t}$. This also
shows that $\diam_{2t}(\Core)$
is bounded by a different constant
times $\log\log{t}$. We now give the details.

\begin{proposition}[The diameter of the inner core]
\label{prop-inner}
Fix $m\geq 2$ and $\delta\in (-m,0)$. Then \whp,
\eq
\label{diam2tinner}
\diam_{2t}(\Inner)\le \frac{2(\tau-1)}{3-\tau}+6.
\en
\end{proposition}

\proof
We first introduce the important notion
of a $t$-connector between a vertex $i\in [t]$ and a set of vertices $A\subseteq[t]$.
This notion will play a crucial role throughout the proof.
We say that the vertex $j\in [2t]\setminus [t]$ is a
{\it $t$-connector} between $i$ and $A$ if
one of the {\it first two} edges incident to $j$ connects to
$i$ and the other of the {\it first two} edges incident to $j$ connects
to a vertex in $A$. Thus, when there exists a $t$-connector between
$i$ and $A$, the distance between $i$ and $A$ in $G_{m,\delta}(2t)$
is at most 2.

We continue the analysis by first considering model (c). We note that for a set of vertices
$A$ and a vertex $i$ with degree at time $t$ equal to
$D_i(t)$, we have that, conditionally on $G_{m,\delta}(t)$, the probability that $j\in [2t]\setminus [t]$ is
a $t$-connector for $i$ and $A$ is at least
    \eq
    \label{tconnprob}
    \frac{(D_{\sss A}(t)+\delta|A|)(D_i(t)+\delta)}{[2t(2m+\delta)]^2}\geq
    \frac{\eta D_{\sss A}(t)D_i(t)}{t^2},
   \en
where in the inequality, we use \ch{that $D_i(t)\geq m$}, and
we let $\eta=(m+\delta)^2/(2m(2m+\delta))^2>0$, while, for any $A\subseteq [t]$, we write
    \eq
    D_{\sss A}(t)=\sum_{i\in A} D_i(t).
    \en
Note that for fixed $j\in [2t]\setminus [t]$ the lower bound
\eqref{tconnprob} holds independently of the fact whether the other vertices
are $t$-connectors or not.

We now give a coupling proof which shows that a subset of size
$n_t=\lfloor \sqrt{t} \rfloor$ of the set $\Inner$ has, \whp, a
bounded \ch{diameter. Lemma} \ref{totdegree} in the appendix
shows that, {\bf whp}, $\Inner$ contains at least $\sqrt{t}$ vertices.
Denote the first $\lfloor \sqrt{t}\rfloor$ vertices of $\Inner$ by
$I$. For each pair $i_1,i_2 \in I$ and each $j\in [2t]\setminus [t]$,
the probability that $j$ is a $t$-connector for $i_1,i_2$ is,
by \eqref{tconnprob}, at least
    \eq
    \label{multprob}
    \frac{\eta u_1^2}{t^2}=\frac{\eta
    t^{\frac1{\tau-1}}}{t^2\log{t}}\geq \frac{
    t^{\frac1{\tau-1}-2}}{\log^2{t}}= q_t,
    \en
independently of the
fact whether the other vertices are $t$-connectors or not. In the
coupling we intend to compare the set $I$ and all pairs
of vertices of the set $I$, which are $t$-connected by some
$j\in [2t]\setminus [t]$ with a so-called {\it multinomial} random graph
$H_{n_t}$. The graph $H_{n_t}$ has $n_t$ vertices and we identify
the $e_t=n_t(n_t-1)/2\sim t/2$ pairs of vertices, which we number
from $1$ to $e_t$ in an arbitrary order, with $e_t$ cells of a
multinomial experiment with $t$ trials and probabilities given by
    \eq
    \label{multinomial}
    p_k=q_t,\quad 1\le k\le e_t,
    \qquad p_0=1-e_tq_t.
    \en
We can represent the $t$ trials by independent random vectors
$N_1,N_2,\ldots,N_t$, where
    \eq
    N_j=(N_{j,1},N_{j,2},\ldots,N_{j,e_t}), \quad 1 \le j\le t,
    \en
with distribution
    \eq
    \prob(N_j=1_{i})=q_t, \quad \prob(N_j=0)=1-e_tq_t,
    \en
where $1_i$ is the $i^{\rm th}$ unit vector of length $e_t$, and $0$ the null vector.
If cell $k$ of the multinomial experiment is not empty,
i.e., if $\sum_{j=1}^t N_{j,k}>0$, then we draw the edge with
number $k$ in the graph $H_{n_t}$, if the cell is empty then
this edge is left out. Note that cell $0$ is just an
overflow cell, which counts the number
of trials that not resulted in one of the cells $1,2,\ldots,e_t$.

By the statement in \eqref{multprob} the distance in $G_{m,\delta}(2t)$ between
any two vertices in $I$ is at most two times the distance between
the corresponding vertices in $H_{n_t}$.
In Lemma \ref{lemmamulti} of the appendix we will show that
the diameter of $H_{n_t}$ is at most the
diameter of a {\it uniform} Erd\H{o}s-R\'enyi graph
$G(n_t,m_t)$, with $n_t$ vertices and $m_t$ edges, where
    \eq
    m_t=\frac12 e_t\Big(1-(1-q_t)^t\Big).
    \en
From \cite[Section 1.4]{JanLucRuc00} we conclude that the above mentioned uniform
Erd\H{o}s-R\'enyi graph $G(n_t,m_t)$ is asymptotically equivalent with the
classical {\it binomial} Erd\H{o}s-R\'enyi graph $G(n_t,\lambda_t)$, where
the edge probability $\lambda_t$ is defined by
    \eq
    \label{edgeprobER}
    \lambda_t=\frac12 \Big(1-(1-q_t)^t\Big)
    \sim \frac{ t^{\frac{1}{\tau-1}-1}}{2\log^2{t}}.
    \en
Next, we show that $\diam(G(n_t, \lambda_t))$ is, \whp,
bounded by $\frac{\tau-1}{3-\tau}+1$. \ch{For this we use the results in
\cite[Corollaries 10.11 and 10.12]{Boll01}, which give sharp bounds on the
diameter of an Erd\H{o}s-R\'enyi random graph. Indeed, this
results imply that if $p^2 n - 2 \log n \rightarrow \infty$ and $n^2(1-p)\rightarrow \infty$, then $\diam(G(n,p))=2$, \whp, while, for $d\geq3$, if $(\log n)/d - 3 \log \log n \rightarrow \infty$ and $p^d n^{d-1}-2\log{n}\rightarrow \infty$,
while $p^{d-1}n^{d-2}-2\log{n}\rightarrow -\infty$, then
$\diam(G(n,p))=d$, \whp. In our case, $n=n_t=\lfloor t^{1/2}\rfloor$ and
$p=\lambda_t$, which implies that, \whp,
$\diam(G(n,p))= \lfloor \frac{\tau-1}{3-\tau}+1\rfloor$.}
We therefore obtain that the diameter of $I$ in $G_{m,\delta}(2t)$ is, \whp, bounded by
    \eq
    \diam_{2t} (I) \leq\frac{2(\tau-1)}{3-\tau}+2.
    \en

We finally show that for any $i\in \Inner\setminus I$, the
probability that there does not exist a $t$-connector connecting
$i$ and $I$ is small. Indeed, since
$D_{\sss I}(t)\geq \sqrt{t} u_1$ and
$D_i(t)\geq u_1$, the mentioned probability is bounded above by
    \eq
    \Big( 1- \frac{\eta D_{\sss I}(t)D_i(t)}{t^2}\Big)^t
    \le
    \exp\left\{
    -\frac{\eta D_{\sss I}(t)D_i(t)}{t}
    \right\}
    \leq
    \exp\left\{
    -\frac{\eta u_1^{2}}{\sqrt{t}}
    \right\}
    \leq
    \exp\left\{
    -\frac{\eta t^{\frac1{\tau-1}-\frac12}}{\log{t}}
    \right\} =o(t^{-1}),
    \en
for $\tau<3$. Thus, {\bf whp}, such a vertex $i$ does not exist.
This proves that {\bf whp} the distance between any vertex
$i\in\Inner \setminus
I$ and $I$ is bounded by 2, and, together with the above bound on
$\diam_{2t}(I)$ we thus obtain \eqref{diam2tinner}.
\qed

\begin{proposition}[Distance from the core to the inner core]
\label{prop-innerouter} Fix $m\geq 2$ and $\delta\in (-m,0)$.
With high probability, the inner core
$\Inner$ can be reached from any vertex in the core $\Core$
using no more than $\frac{2\log\log{t}}{|\log{(\tau-2)}|}$ edges in $G_{m,\delta}(2t)$. More
precisely, {\bf whp},
    \eq
    \max_{i\in \Core} \min_{j\in \Inner} \dist_{\sss G_{m,\delta}(2t)}(i,j)\leq
    \frac{2 \log\log{t}}{|\log{(\tau-2)}|}.
    \en
\end{proposition}

\proof For $k\geq 1$, we define
    \eq
        \Ncal{k}=\{i\in[t]: D_i(t)\geq u_k\},
    \en
with $u_1$ defined in \eqref{Inner-def}, and where
we define $u_k,$ for $k\ge 2$, recursively, so that
for any vertex $i\in [t]$ with degree at least $u_k$,
the probability that there is no $t$-connector for the vertex $i$
and the set $\Ncal{k-1}$, conditionally on $G_{m,\delta}(t)$, is tiny. According to
\eqref{tconnprob} and \eqref{lowerbddegrees} in the appendix,
this probability is at
most
    \eq
    \label{fourtwenty}
    \left(
    1-\frac{\eta \ch{D}_{\sss \Ncal{k-1}} D_i(t)}{t^2}
    \right)^t
    \le
     \exp\Big\{-\frac{\eta B t\big(u_{k-1}\big)^{2-\tau}u_k}{t}\Big\}=o(t^{-2}),
    \en
for some $B>0$, when we define
    \eq
        \label{Wk-def}
        \Wcal{k}=D\log{t} \big(\Wcal{k-1}\big)^{\tau-2},
    \en
with $D$ exceeding $2(\eta B)^{-1}$ \ch{ and $t$ is sufficiently large so that $\Wcal{k}\leq \Wcal{1}$}. The following lemma
identifies $\Wcal{k}$:

\begin{lemma}[Identification of $\Wcal{k}$]
\label{lem-powersWk}
For each $k\in \N$,
    \eq
    \Wcal{k}=D^{a_k} (\log{t})^{b_k} t^{c_k},
    \en
where
\eq\label{gerard10}
    a_k=\frac{1-(\tau-2)^{k-1}}{3-\tau},\qquad
    b_k=\frac{1-(\tau-2)^{k-1}}{3-\tau}-\frac{1}{2}(\tau-2)^{k-1},
    \qquad
    c_k=\frac{(\tau-2)^{k-1}}{2(\tau-1)}.
\en
\end{lemma}

\proof We leave the straightforward induction proof to the reader.
\qed
\vskip0.5truecm

Then, the key step in the proof of Proposition
\ref{prop-innerouter} is the following lemma:

\begin{lemma}[Connectivity between $\Ncal{k-1}$ and $\Ncal{k}$]
\label{lem-connNs} Fix $m \geq 2$ and $\delta\in (-m,0)$. Then,
uniformly in $k$, the probability that there
exists an $i\in \Ncal{k}$ that is not at distance
at most two from
$\Ncal{k-1}$ in $G_{m,\delta}(2t)$ is $o(t^{-1})$.
\end{lemma}

\proof It follows from \ch{\eqref{fourtwenty}} that the probability in
the statement is by Boole's inequality bounded by
    \begin{align}
    t \exp\Big(
    -\frac{\eta B t[\Wcal{k-1}]^{2-\tau}\Wcal{k}}{t}\Big)
    &=t \cdot o(t^{-2})=o(t^{-1})\ch{.}
    \end{align}
\qed
\vskip0.5cm

\noindent We now complete the proof of Proposition
\ref{prop-innerouter}. Fix
    \eq
    \label{k*-def}
    k^*= \Big\lfloor \frac{\log\log{t}}{|\log{(\tau-2)}|}\Big \rfloor.
    \en
As a
result of Lemma \ref{lem-connNs}, we have that the distance
between $\Ncal{k^*}$ and $\Inner= \Ncal{1}$ is at most $2k^*$. Therefore, Proposition~\ref{prop-innerouter} follows when we can show that
    \eq
    \Core= \{i: D_i(t) \geq (\log{t})^\sigma\} \subseteq
    \Ncal{k^*}=\{i: D_i(t) \geq \Wcal{k^*}\},
    \en
so that it suffices to prove that $(\log{t})^{\sigma}\geq \Wcal{k^*}$,
for any $\sigma>1/(3-\tau)$.
This follows trivially for $t$ large from the explicit
representation of $u_{k^*}$
given by Lemma \ref{lem-powersWk}.
\qed \vskip0.5cm

\noindent {\it Proof of Theorem \ref{prop-core}.} We note that
\whp~
\eq
\diam_{2t}(\Core)\leq  \frac{2(\tau-1)}{3-\tau}+6+\ch{4}k^*,
\en
where $k^*$ is given in \eqref{k*-def}, and where we have made use of
\ch{Propositions \ref{prop-inner} and \ref{prop-innerouter}}. This proves Theorem \ref{prop-core}.
\qed

\subsection{Connecting the periphery to the core}
\label{sec-percore}
In this section, we extend the results
of the previous section and, in particular,
study the distance between the vertices not in the core
$\Core$ and the core. The main result is
the following theorem:

\begin{theorem}[Connecting the periphery to the core]
\label{prop-periphery}
Fix $m\geq 2$ and $\delta\in (-m,0)$. For every
$\sigma>1/(3-\tau)$, {\bf whp},
the maximal distance between
any vertex and  $\Core$ in $G_{m,\delta}(2t)$ is bounded from above by
$2\sigma\log\log{t}/\log{m}$.
\end{theorem}

\noindent Together with Theorem \ref{prop-core}, Theorem \ref{prop-periphery}
proves the main result in Theorem \ref{thm-diametersec}.

The proof of Theorem \ref{prop-periphery} consists
of two key steps. The first key step in Proposition \ref{prop-halfway}
states that the distance between any vertex in $[t]$ and the core $\Core$ is bounded by
a constant times $\log\log{t}$.
The second key step in Proposition \ref{prop-comple} shows that
the distance between any vertex in $[2t]\setminus [t]$ and
$[t]$ is bounded by another constant times $\log\log{t}$.

\begin{proposition}[Connecting half of the periphery to the core]
\label{prop-halfway}
Fix $m\geq 2$ and $\delta\in (-m,0)$. For every $\sigma>1/(3-\tau)$,
{\bf whp}, the distance between
any vertex in $[t]$ and the core $\Core$ in $G_{m,\delta}(2t)$ is bounded from above by
$\sigma\log\log{t}/\log{m}$.
\end{proposition}

\proof
We start from a vertex $i\in [t]$ and will show that
the probability that the distance between $i$ and $\Core$
is at least $\sigma\log\log{t}/\log{m}$ is $o(t^{-1})$.
This proves the claim. For this, we explore the neighborhood of
$i$ as follows. From $i$, we connect its $m\geq 2$ edges.
Then, successively, we connect the $m$
edges from each of the at most $m$ vertices that $i$ has connected to and
have not yet been explored. We continue in the same fashion.
We call the arising process when we
have explored up to distance $k$ from the initial
vertex $i$ the {\it $k$-exploration tree} of vertex $i$.

When we never connect two edges to the same vertex,
then the number of vertices we can reach within $k$ steps is {\it precisely}
equal to $m^{k}$. We call an event where an edge connects to
a vertex which already was in the exploration tree a {\it collision}.
When $k$ increases, the probability
of a collision increases. However, the probability that
there exists a vertex for which {\it more than \ch{$l$}}
collisions occur in its $k$-exploration tree\ch{, where $\l\geq1$,} before it hits the core
is small, as we prove now:

\begin{lemma} [A bound on the probability of multiple collisions]
\label{lem-collision}
Fix $m\geq 2$ and $\delta\in (-m,0)$. Fix $l\geq 1$, $b\in(0,1]$
and take $k\le \sigma\log\log{t}/\log{m}$. Then, for every vertex $i\in [t]$,
the probability that its $k$-exploration tree has at
least $l$ collisions before it hits $\Core\cup [t^b]$
is bounded above by
    \eq
    m^l(\log{t})^{2 \sigma l} \Big/  t^{bl}.
    \en
\end{lemma}

\proof Take $i\in [t] \setminus [t^b]$ and consider its
$k$-exploration tree $\Tcal{k}{i}$.
Since we add edges after time $t^b$ the denominator in \eqref{growthrulePA}-\eqref{growthrulePAc}
is at least $t^b$. Moreover, before hitting the core, any vertex in the $k$-exploration tree
has degree at most $(\log t)^{\sigma}$. Hence, for $l=1$, the probability mentioned in the
statement of the lemma is at most
    \eq
    \label{gerard11}
    \sum_{v\in \Tcal{k}{i}} \frac{D_v(t)+\delta}{t^b}\leq \sum_{v\in \Tcal{k}{i}} \frac{(\log t)^{\sigma}}{t^b}
    \leq \frac{m^{k+1}(\log{t})^{\sigma}}{t^b},
    \en
where the bound follows from $\delta<0$ and $\#\{v\in\Tcal{k}{i}\}\le m^{k+1}$.
For general $l$ this upper bound becomes:
    \eq
    \left(\frac{m^{k+1}(\log{t})^{\sigma}}{t^b}\right)^l.
    \en
When $k=\sigma\log\log{t}/\log{m}$,
we have that $m^{kl}= (\log{t})^{\sigma l}$.
Therefore, the claim in Lemma \ref{lem-collision}
holds.
\qed
\vskip0.5cm

\noindent
We next prove that there exists a $b>0$ such that, \whp, $[t^b]$
is a subset of the core. Note that in this lemma the conditions $m\ge 2$ or $\delta\in (-m,0)$
are not necessary.

\begin{lemma}[Early vertices have large degrees \whp]
\label{lem-earlypartofcore}
Fix $m\geq 1$. There exists a $b>0$ such that, \ch{for every $\sigma>1/(3-\tau)$,
\whp, $\min_{j\leq t^b} D_j(t)\geq (\log{t})^{\sigma}$.} As a result, \whp,
$[t^b]\subseteq \Core$.
\end{lemma}

We defer the proof of Lemma \ref{lem-earlypartofcore} to
Section \ref{sec-prfslems} of the appendix. Now we are
ready to complete the proof of Proposition \ref{prop-halfway}:
\vskip0.2cm

\noindent
{\it Proof of Proposition \ref{prop-halfway}.}
By combining Lemmas \ref{lem-collision}
and \ref{lem-earlypartofcore}, the probability that there exists an
$i\in [t]$ for which the exploration tree $\Tcal{k}{i}$ has at least
$l$ collisions before hitting the core is $o(1)$, whenever $l>1/b$, since,
by Boole's inequality, it is bounded by
    \eq
    m^l \sum_{i=1}^t (\log{t})^{2\sigma l} \Big/ t^{bl}=
    m^l(\log{t})^{2\sigma l} t^{-bl+1}=o(1)\ch{.}
    \en
\ch{When} the $k$-exploration tree hits the core, then we are done.
When the $k$-exploration tree from
a vertex $i$ does not hit the core,
but has less than $l$ collisions, then
there are at least $m^{k-l}$ vertices in $k$-exploration tree.
\ch{Indeed, when we have at most $l$ collisions, the size of the $k$-exploration tree is minimal when all edges of the root connect to the same vertex $v_1$, all edges of $v_1$ connect to the same vertex $v_2$, etc. Iterating this at most $l$ levels deep yields a tree with at least $m^{k-l}$ vertices.}

When $k=\sigma\log\log{t}/\log{m}-2$,
$m^{k-l}\geq (\log{t})^{\sigma+o(1)}$.
The total \ch{degree} of the core is, by \eqref{lowerbddegrees} in the appendix,
at least
    \eq\label{eqrev1}
    \sum_{i\in \Core} \ch{D_i(t)}\geq B t (\log{t})^{-(\tau-2)\sigma},
    \en
for some $B>0$.
The probability that there does not exist a $t$-connector
between the $k$-exploration tree and the core is, by
\eqref{tconnprob} and \ch{\eqref{eqrev1},} bounded above by
    \eq
    \exp\left\{-\frac{\eta B t (\log{t})^{-(\tau-2)\sigma}
    (\log{t})^{\sigma+o(1)}}{t}\right\}
    =o(t^{-1}),
    \en
since $\sigma > 1/(3-\tau)$. This completes the proof.
\qed

\begin{proposition}[Connecting the remaining periphery]
\label{prop-comple}
Fix $m\geq 2$ and $\delta \in(-m,0)$.
For every $\sigma >1/(3-\tau)$, {\bf whp}, the maximal distance between
any vertex and $[t]$ in $G_{m,\delta}(2t)$ is bounded from above by
$\sigma\log\log{t}/\log{m}$.
\end{proposition}

\proof Take $k=\sigma\log\log{t}/\log{m}-1$, and
$j\in[2t]\setminus [t]$ with distance larger than $k$ to the set of vertices
$[t]$.
We now apply  Lemma \ref{lem-collision}
with $t$ replaced by $2t$ and letting $l=2$ and \ch{$b=b_t\in(0,1)$} such that $(2t)^b=t$, to conclude that
with probability exceeding $1-o(1)$, the $k$-exploration tree of $j$ has at most $1$ collision
before it hits ${\rm Core}_{2t} \cup [t]$.
We can hence conclude that with probability exceeding $1-o(1)$, there are at least
$m_k=(m-1)m^{k-1}$ vertices in $[2t] \setminus [t] $ at distance precisely equal to $k$ from our starting vertex $j$.
Denote these vertices by $i_1, \ldots, i_{m_k}$.
We consider case (c), the proof for (a) and (b) is similar.
Note that, uniformly in $s\in [2t]\setminus [t]$,
    \eq
    \frac{\sum_{i=1}^t (D_i(s)+\delta)}{(2m+\delta)s}\geq \frac 12.
    \en
Hence,
    \eq
    \prob\Big(\nexists l\in [m_k] \text{ such that }
    {\rm dist}_{G_{m,\delta}(2t)}(i_l, {\rm Core}_{2t} \cup [t]) \ch{\leq} 1 \Big)\\
    \leq 2^{-m_k}=o(t^{-1}),
    \en
since $m_k=\frac{m-1}{m^2}(\log{t})^{\sigma}$, with $\sigma>1/(3-\tau)>1$.
Therefore, any vertex $j\in [2t]\setminus [t]$ is, {\bf whp},
within distance $k+1$ from ${\rm Core}_{2t} \cup [t]$. Proposition~\ref{latevertsmalldeg}
shows that, {\bf whp} the set
${\rm Core}_{2t}\subseteq [t]$, so that, {\bf whp}, ${\rm Core}_{2t}\cup [t]=[t]$
and the proposition follows.
\qed
\vskip0.5cm

\noindent
{\it Proof of Theorem \ref{prop-periphery}.}
Proposition \ref{prop-comple} states that {\bf whp} every vertex in $G_{m,\delta}(2t)$ is
within distance $\sigma\log\log{t}/\log{m}$
of $[t]$ and Proposition \ref{prop-halfway} states that
{\bf whp} every vertex in $[t]$ is at most distance $\sigma\log\log{t}/\log{m}$
from the core $\Core$. This shows that every vertex in $G_{m,\delta}(2t)$
is {\bf whp} within distance $2\sigma\log\log{t}/\log{m}$ from the core.
\qed
\vskip0.5cm

\noindent
{\it Proof of Theorem \ref{thm-diametersec}.}
Theorem \ref{prop-periphery} states that every vertex in $G_{m,\delta}(2t)$ is
within distance $\frac{2\sigma\log\log{t}}{\log{m}}$
of the core $\Core$. Theorem \ref{prop-core}
states that the diameter of the core is at most
$\frac{4 \log\log{t}}{|\log{(\tau-2)}|}(1+o(1))$, so that the diameter
of $G_{m,\delta}(2t)$ is at most $C_{\sG}\log\log{t}$, where $C_{\sG}$
is given in \refeq{CG-def}, because we can choose any $\sigma > 1/(3-\tau)$.
This completes the proof of Theorem \ref{thm-diametersec}.
\qed

\section{A $\log\log{t}$ lower bound on the diameter: Proof of Theorem \ref{thmlowerb}}
\label{lowerbound}
We will again prove this theorem for time $2t$ rather than time $t$.
To show that the diameter of the graph is, {\bf whp}, at least $k$, we will study, at time $2t$,
the $k$-exploration trees $\mathcal{T}_i^{\sss(k)}$ of vertices $i\in [2t]\backslash[t]$ as defined above.
We shall call the tree $\mathcal{T}_i^{\sss(k)}$ \emph{proper} if the following conditions hold:
\begin{itemize}
\item The $k$-exploration tree has no collisions;
\item All vertices of $\mathcal{T}_i^{\sss(k)}$ are in $[2t]\backslash[t]$;
\item No other vertex connects to a vertex in $\mathcal{T}_i^{\sss(k)}$.
\end{itemize}

When such a tree exists in $G_{m,\delta}(2t)$ for a certain vertex $i$ then we know that the diameter is at least $k$,
since the distance between the root of the tree $i$ and the vertices at depth $k$ is exactly $k$; there cannot be a shorter route.

To prove that a proper $k$-exploration tree exists in $G_{m,\delta}(2t)$, we will use a second moment method.
Let $\mathfrak{T}_m^k(2t)$ be the set of all possible $k$-exploration trees that can exist in $G_{m,\delta}(2t)$
and satisfy the first two conditions.
Note that the order in which the edges are added matters: if two edges are added in a different order,
then the arising exploration tree will be considered a different tree.
Let $Z_{m,\delta}^{\sss(k)}(2t)$ be the number of proper $k$-exploration trees in $G_{m,\delta}(2t)$, i.e.,
\begin{equation}
Z_{m,\delta}^{\sss(k)}(2t) = \sum_{\mathcal{T}\in\mathfrak{T}_m^k(2t)} I[ \mathcal{T} \subseteq G_{m,\delta}(2t) \text{ and }
\mathcal{T} \text{ is proper} ].
\end{equation}
Here the event that all edges of $\mathcal{T}$ have been formed in $G_{m,\delta}(2t)$ is denoted by $\mathcal{T} \subseteq G_{m,\delta}(2t)$.

In Section~\ref{secfirst} we will investigate the first moment of $Z_{m,\delta}^{\sss(k)}(2t)$ and prove the following:
\begin{proposition}[Expected number of proper trees tends to infinity] \label{propfirst}
Fix $m\geq2$ and $\delta>-m$. Let $k=\frac{\varepsilon}{\log m} \log\log t$, with $0<\varepsilon<1$. Then
\begin{equation}
\lim_{t\rightarrow\infty} \mathbb{E}\left[Z_{m,\delta}^{\sss(k)}(2t)\right] = \infty.
\end{equation}
\end{proposition}

The variance of $Z_{m,\delta}^{\sss(k)}(2t)$ will be the subject of Section~\ref{secsecond}, where we will prove the following:
\begin{proposition}[Concentration of the number of proper trees] \label{propsecond}
Fix $m\geq2$, $\delta >-m$ and \ch{let} $0 \leq k \leq \frac{\log\log t}{\log m}$.
Then there exists a constant $c_{m,\delta}>0$, such that, for $t$ sufficiently large,
\begin{equation}
\Var\left(Z_{m,\delta}^{\sss(k)}(2t)\right) \leq c_{m,\delta} \frac{(\log t)^2}{t} \mathbb{E}\left[Z_{m,\delta}^{\sss(k)}(2t)\right]^2
+\mathbb{E}\left[Z_{m,\delta}^{\sss(k)}(2t)\right].
\end{equation}
\end{proposition}
We use these two propositions to prove Theorem \ref{thmlowerb}:

\begin{proof}[Proof of Theorem~\ref{thmlowerb}]
We first use the Chebychev inequality to obtain that
\begin{align}\label{eqcheb}
\mathbb{P}\left(\text{diam}(G_{m,\delta}(2t)) < k\right)&\leq \mathbb{P}\left(Z_{m,\delta}^{\sss(k)}(2t) = 0\right)
\leq\frac{\Var\left(Z_{m,\delta}^{\sss(k)}(2t)\right)}{\mathbb{E}\left[Z_{m,\delta}^{\sss(k)}(2t)\right]^2}.
\end{align}
By Proposition~\ref{propsecond}, the right-hand side of \eqref{eqcheb} is, for some constant $c_{m,\delta}>0$, at most
\begin{equation}
c_{m,\delta} \frac{(\log t)^2}{t} +\frac{1}{\mathbb{E}\left[Z_{m,\delta}^{\sss(k)}(2t)\right]}=o(1),\label{eqvarb}
\end{equation}
by Proposition~\ref{propfirst}.
\end{proof}

\subsection{The first moment of the number of proper trees}
\label{secfirst}
Let $\mathcal{B}_\mathcal{T}$ denote the event that no vertex outside a tree $\mathcal{T}$ connects to a vertex in this tree.
We can then write that the expected number of proper $k$-exploration trees in $G_{m,\delta}(2t)$ equals
\begin{align}
\mathbb{E}\left[Z_{m,\delta}^{\sss(k)}(2t)\right]
&= \sum_{\mathcal{T}\in\mathfrak{T}_m^k(2t)}
\mathbb{P}\Big( \mathcal{T} \subseteq G_{m,\delta}(2t) \text{ and } \mathcal{T} \text{ is proper} \Big) \nonumber\\
&= \sum_{\mathcal{T}\in\mathfrak{T}_m^k(2t)}
\mathbb{P}\Big(\mathcal{T} \text{ is proper}| \mathcal{T} \subseteq G_{m,\delta}(2t)\Big)\mathbb{P}\Big( \mathcal{T} \subseteq G_{m,\delta}(2t) \Big)  \nonumber\\
&= \sum_{\mathcal{T}\in\mathfrak{T}_m^k(2t)}
\mathbb{P}\Big(\mathcal{B}_\mathcal{T}| \mathcal{T} \subseteq G_{m,\delta}(2t)\Big) \cdot\mathbb{P}\Big( \mathcal{T} \subseteq G_{m,\delta}(2t) \Big).
\end{align}

We will first give a lower bound on the probability that a given $k$-exploration tree exists in the graph at time $2t$. For convenience we will write $a_{m,\delta}=\frac{m+\delta}{3(2m+\delta)}$.
\begin{lemma}[Lower bound on existence probability]\label{lemptree}
Fix $m\geq2$, $\delta > -m$ and $k\geq0$. Given a proper $k$-exploration tree $\mathcal{T}\in\mathfrak{T}_m^k(2t)$, then, for $t$ sufficiently large,
\begin{equation}
\mathbb{P}\Big(\mathcal{T}\subseteq G_{m,\delta}(2t)\Big) \geq \left(\frac{a_{m,\delta}}{t}\right)^{m^{\sss(k)}-1},
\end{equation}
where $m^{\sss(k)}=\frac{m^{k+1}-1}{m-1}$.
\end{lemma}
\begin{proof}
Since every vertex is added before time $2t$, the denominator in \eqref{growthrulePA}--\eqref{growthrulePAc}  is at most
$3t(2m+\delta)$. The degree of all vertices already in the graph is at least $m$, so the probability that a certain given edge is formed is at least
\begin{equation}
\frac{m+\delta}{3t(2m+\delta)}=\frac{a_{m,\delta}}{t}.
\end{equation}
Since exactly $m^{\sss(k)}-1$ edges have to be formed to form the given tree $\mathcal{T}$, we have that
\begin{equation}
\mathbb{P}\Big(\mathcal{T}\subseteq G_{m,\delta}(2t)\Big) \geq\left(\frac{a_{m,\delta}}{t}\right)^{m^{\sss(k)}-1}.
\end{equation}
\end{proof}

We will now give a lower bound on the probability that no other vertex connects to a given tree.
\begin{lemma}[No other vertex connects to $\mathcal{T}$]\label{lempproper}
Fix $m\geq2$, $\delta >-m$ and $0\leq k \leq \frac{\log\log t}{\log m}$. Given a proper $k$-exploration tree $\mathcal{T}\in\mathfrak{T}_m^k(2t)$, then, for $t$ sufficiently large, \ch{and writing $m_\delta=m+1+\delta>1$,}
\begin{equation}
\mathbb{P}\Big(\mathcal{B}_\mathcal{T}| \mathcal{T} \subseteq G_{m,\delta}(2t)\Big) \geq \left(1-\frac{m_\delta m^{k+1}}{t}\right)^{mt}.
\end{equation}
\end{lemma}
\ch{
\begin{rem}
In $\mathbb{P}\Big(\mathcal{B}_\mathcal{T}| \mathcal{T} \subseteq G_{m,\delta}(2t)\Big)$, $\mathcal{B}_\mathcal{T}$ makes a claim about edges {\em not} in $\mathcal{T}$, while the event $\mathcal{T} \subseteq G_{m,\delta}(2t)$ states that all edges in $\mathcal{T}$ are formed in our random graph process. Thus conditioning on $\mathcal{T} \subseteq G_{m,\delta}(2t)$ gives information only about {\em inside} edges.
\end{rem}
}
\begin{proof}
First note that for $k \leq \frac{\log\log t}{\log m}$ and $t$ sufficiently large, $m_\delta m^{k+1} \leq m_\delta m\log t \leq t$. So $0\leq 1-\frac{m_\delta m^{k+1}}{t}\leq1$. Further note that vertices in $[t]$ cannot connect to a vertex in $\mathcal{T}$, since $\mathcal{T}\subseteq[2t]\backslash[t]$. In the remainder of the proof we will refer to \emph{outside} edges as those edges that do not belong to $\mathcal{T}$, of which there are exactly $mt-(m^{\sss(k)}-1)$ added after time $t$. For $A$ a set of vertices, let $\mathcal{E}_n(A)$ denote the event that the $n$-th outside edge added after time $t$ connects to a vertex in $A$ and let $\overline{\mathcal{E}}_n(A)$ be the complement of $\mathcal{E}_n(A)$. We use induction on the number of outside edges that did not connect to the tree $\mathcal{T}$, i.e., we show that:
\begin{equation}
\mathbb{P}\left(\bigcap_{i=1}^n \overline{\mathcal{E}}_i(\mathcal{T})\Big| \mathcal{T} \subseteq G_{m,\delta}(2t)\right) \geq \left(1-\frac{m_\delta m^{k+1}}{t}\right)^{n},
\end{equation}
by induction on $n=0,\ldots,mt-(m^{\sss(k)}-1)$. For $n=0$ the above holds, because both sides equal $1$. Now assume that the above holds for $0\leq n <mt-(m^{\sss(k)}-1)$, then
\begin{align}
\mathbb{P}\Bigg(&\bigcap_{i=1}^{n+1} \overline{\mathcal{E}}_i(\mathcal{T})\Big| \mathcal{T} \subseteq G_{m,\delta}(2t)\Bigg)\nonumber\\
&=\mathbb{P}\left(\overline{\mathcal{E}}_{n+1}(\mathcal{T})\Big|\bigcap_{i=1}^{n} \overline{\mathcal{E}}_i(\mathcal{T})\cap \left\{\mathcal{T} \subseteq G_{m,\delta}(2t)\right\}\right) \mathbb{P}\left( \bigcap_{i=1}^{n} \overline{\mathcal{E}}_i(\mathcal{T})\Big| \mathcal{T} \subseteq G_{m,\delta}(2t)\right)\nonumber\\
&\geq\left(1-\mathbb{P}\left( \mathcal{E}_{n+1}(\mathcal{T})\Big|\bigcap_{i=1}^{n} \overline{\mathcal{E}}_i(\mathcal{T})\cap \left\{\mathcal{T} \subseteq G_{m,\delta}(2t)\right\}\right)\right)\cdot \left(1-\frac{m_\delta m^{k+1}}{t}\right)^{n}.\label{eqinduc}
\end{align}
Since it is known that at the time that the $(n+1)$-st outside edge after time $t$ is added,
no other outside edge has connected to a vertex in the tree, we know that the degree of all vertices in the tree at that moment is at most $m+1$.
Further, since this edge is added after time $t$, the denominator of
\eqref{growthrulePA}--\eqref{growthrulePAc}
 will be at least $t$. Thus, the right-hand side of \eqref{eqinduc} is at least
\begin{align}
\left(1-\sum_{i\in\mathcal{T}}\frac{m+1+\delta}{t}\right)\cdot \left(1-\frac{m_\delta m^{k+1}}{t}\right)^{n}
&\geq\left(1-\frac{m_\delta m^{k+1}}{t}\right)\cdot\left(1-\frac{m_\delta m^{k+1}}{t}\right)^{n}\nonumber\\
&=\left(1-\frac{m_\delta m^{k+1}}{t}\right)^{n+1},
\end{align}
where the inequality holds because there are less than $m^{k+1}$ vertices in the tree. Applying the above to $n=mt-(m^{\sss(k)}-1)$, we obtain that
\begin{equation}
\mathbb{P}\Big(\mathcal{B}_\mathcal{T}| \mathcal{T} \subseteq G_{m,\delta}(2t)\Big) \geq \left(1-\frac{m_\delta m^{k+1}}{t}\right)^{mt-(m^{\sss(k)}-1)}\geq \left(1-\frac{m_\delta m^{k+1}}{t}\right)^{mt}.
\end{equation}
\end{proof}

We finally give a lower bound on the number of possible proper $k$-exploration trees that can be formed. It should be noted that when a vertex $i$ connects to a vertex $j$, we will always have that $i>j$. So when exploring a vertex $i$ in the exploration tree, all $m$ vertices this vertex connects to have a smaller label than $i$.

\begin{lemma}[Number of proper trees] \label{lemntrees}
Fix $m\geq2$ and $0\leq k\leq \frac{\log\log t}{\log m}$. Then, for $t$ sufficiently large, the number of possible proper $k$-exploration trees at time $2t$ is at least $\left(t/m^{k+1}\right)^{m^{\sss(k)}}$,
where we recall that $m^{\sss(k)}=\frac{m^{k+1}-1}{m-1}$.
\end{lemma}

\begin{proof}
For $t$ sufficiently large and $ k\leq \frac{\log\log t}{\log m}$, $m^{k+1}\leq m\log t \leq t$. Since the $k$-exploration tree of a vertex $i$ has to be proper, there are no collisions, so the number of vertices in the tree equals
\begin{equation}
\#\{v\in\mathcal{T}_i^{\sss(k)}\}=m^{\sss(k)}.
\end{equation}
For any subset $X \subseteq [2t]\backslash[t]$ with $\#\{v\in X\}=m^{\sss(k)}$ there exists at least one possible proper $k$-exploration tree. To see this, first order the vertex labels in descending order. Let the first vertex, i.e. the vertex with the largest label, be the root of the tree. Then let the next $m$ vertices be the vertices at distance 1 from the root, the next $m^2$ vertices be the vertices at distance 2 from the root, etcetera, until the last $m^k$ vertices which will be at distance $k$ from the root. This way, all vertices will connect to $m$ vertices with a smaller label, i.e., vertices that were already in the graph when the vertex was added, so this is a possible proper $k$-exploration tree with all vertices in $X$.

The number of subsets of $[2t]\backslash[t]$ of size $m^{\sss(k)}$ is $\binom{t}{m^{\sss(k)}}$ which is at least
\begin{equation}
\left(\frac{t}{m^{\sss(k)}}\right)^{m^{\sss(k)}}\geq\left(\frac{t}{m^{k+1}}\right)^{m^{\sss(k)}},
\end{equation}
where we used that for $1 \leq b\leq a$ we have that $(a-i)b \geq (b-i)a$ for all $0 \leq i < b$, so that
\begin{equation}
\binom{a}{b} = \prod_{i=0}^{b-1}\frac{a-i}{b-i} \geq \left(\frac{a}{b}\right)^b.
\end{equation}
\end{proof}

We can now combine the three bounds above to get a lower bound on the expected number of proper $k$-exploration trees.
\begin{corollary}[Lower bound on expected number of proper trees]
\label{colbound1}
Fix $m\geq2$, $\delta>-m$ and $0\leq k\leq \frac{\log\log t}{\log m}$. Then, for $t$ sufficiently large,
\begin{equation}
\mathbb{E}\left[Z_{m,\delta}^{\sss(k)}(2t)\right] \geq \frac{t}{a_{m,\delta}} \left(\frac{a_{m,\delta}}{m^{k+1}}\right)^{m^{k+1}}  \left(1-\frac{m_\delta m^{k+1}}{t}\right)^{mt}.
\end{equation}
\end{corollary}

\begin{proof}
Using the bounds from Lemmas~\ref{lemptree},~\ref{lempproper} and~\ref{lemntrees} we get that
\begin{align}
\mathbb{E}\left[Z_{m,\delta}^{\sss(k)}(2t)\right] &= \sum_{\mathcal{T}\in\mathfrak{T}_m^k(2t)} \mathbb{P}\Big(\mathcal{B}_\mathcal{T}| \mathcal{T} \subseteq G_{m,\delta}(2t)\Big) \cdot\mathbb{P}\Big( \mathcal{T} \subseteq G_{m,\delta}(2t) \Big)\nonumber\\
&\geq \#\{\mathcal{T}\in\mathfrak{T}_m^k(2t)\} \left(1-\frac{m_\delta m^{k+1}}{t}\right)^{mt} \left(\frac{a_{m,\delta}}{t}\right)^{m^{\sss(k)}-1}\nonumber\\
&\geq \left(\frac{t}{m^{k+1}}\right)^{m^{\sss(k)}}  \left(1-\frac{m_\delta m^{k+1}}{t}\right)^{mt} \left(\frac{a_{m,\delta}}{t}\right)^{m^{\sss(k)}-1}\nonumber\\
&\geq\frac{t}{a_{m,\delta}} \left(\frac{a_{m,\delta}}{m^{k+1}}\right)^{m^{k+1}}  \left(1-\frac{m_\delta m^{k+1}}{t}\right)^{mt}.
\end{align}
\end{proof}
The factor $t$ in the corollary above turns out to be crucial for the remainder of the proof. This factor arises from the fact that there is exactly one edge less in a proper $k$-exploration tree than there are vertices.

We can now show that the expected number of $k$-exploration trees tends to infinity, for $k=\frac{\varepsilon}{\log m} \log\log t$, with $0<\varepsilon<1$.
\begin{proof}[Proof of Proposition~\ref{propfirst}]
First note that for $k=\frac{\varepsilon}{\log m} \log\log t$, with $0<\varepsilon<1$, $m^k=(\log t)^{\varepsilon}$. We can then use Corollary~\ref{colbound1} to get that
\begin{equation}
\lim_{t\rightarrow\infty} \mathbb{E}\left[Z_{m,\delta}^{\sss(k)}(2t)\right] \geq \lim_{t\rightarrow\infty}\frac{t}{a_{m,\delta}}\left(\frac{a_{m,\delta}}{m^{k+1}}\right)^{m^{k+1}}  \left(1-\frac{mm_\delta m^{k+1}}{mt}\right)^{mt} =\infty,
\end{equation}
since
\begin{equation}
\left(\frac{a_{m,\delta}}{m^{k+1}}\right)^{m^{k+1}} = \left(\frac{a_{m,\delta}}{m(\log t)^\varepsilon}\right)^{m(\log t)^\varepsilon},
\quad \text{and} \quad
\left(1-\frac{mm_\delta m^{k+1}}{mt}\right)^{mt} \sim e^{-m^2m_\delta (\log t)^\varepsilon}.
\end{equation}
\end{proof}
It is easy to see that the same argument can be applied to $k=\frac{\log \log t}{\log m}-\frac{\log \log \log t}{\log m}-1$.

\subsection{The second moment of the number of proper trees}
\label{secsecond}
In this section we will investigate the variance of $Z_{m,\delta}^{\sss(k)}(2t)$. To shorten the notation, for a $k$-exploration tree $\mathcal{T}\in\mathfrak{T}_m^k(2t)$, let $F_{\mathcal{T}}$ denote the event that $\mathcal{T} \subseteq G_{m,\delta}(2t) \text{ and } \mathcal{T} \text{ is proper}$. Then, the variance of the number of proper $k$-exploration trees in $G_{m,\delta}(2t)$ is given by
\begin{align}
\Var\left(Z_{m,\delta}^{\sss(k)}(2t)\right) &= \Var\left(\sum_{\mathcal{T}\in\mathfrak{T}_m^k(2t)} I[ \mathcal{T} \subseteq G_{m,\delta}(2t) \text{ and } \mathcal{T} \text{ is proper} ]\right)\nonumber\\
&= \Var\left(\sum_{\mathcal{T}\in\mathfrak{T}_m^k(2t)} I[ F_{\mathcal{T}} ]\right)
= \sum_{\mathcal{T,T'}\in\mathfrak{T}_m^k(2t)}\Cov\left(I[ F_{\mathcal{T}}],I[F_{\mathcal{T'}} ]\right) \nonumber\\
&= \sum_{\substack{\mathcal{T,T'}\in\mathfrak{T}_m^k(2t) \\ \mathcal{T}\neq{\mathcal{T'}}}}\left(\mathbb{P}\left( F_{\mathcal{T}} \cap F_{\mathcal{T'}} \right) -\mathbb{P}\left(F_{\mathcal{T}}\right)\mathbb{P}\left(F_{\mathcal{T'}}\right)\right)
+\sum_{\mathcal{T}\in\mathfrak{T}_m^k(2t)}\mathbb{P}\left(F_{\mathcal{T}}\right)\left(1-\mathbb{P}\left(F_{\mathcal{T}}\right)\right). \end{align}
We start by studying the terms of the first sum in the following lemma.
\begin{lemma}[Weak dependence of tree occurrences] \label{lempminpp}
Fix $m\geq2$, $\delta > -m$ and $0 \leq k \leq \frac{\log\log t}{\log m}$. Let $\mathcal{T,T'}\in\mathfrak{T}_m^k(2t)$ with $\mathcal{T} \neq \mathcal{T'}$. Then, for $t$ sufficiently large,
\begin{equation}\label{eqlemsec}
\mathbb{P}\left( F_{\mathcal{T}} \cap F_{\mathcal{T'}} \right)- \mathbb{P}\left(F_{\mathcal{T}}\right)\mathbb{P}\left(F_{\mathcal{T'}}\right)\leq\left(\left(1+\frac{2m_\delta m \log t}{t}\right)^{2m \log t} -1\right)\mathbb{P}\left(F_{\mathcal{T}}\right)\mathbb{P}\left(F_{\mathcal{T'}}\right).
\end{equation}
\end{lemma}
\begin{proof}
When $\mathcal{T} \cap \mathcal{T'} \neq \emptyset$, at least one edge of one of the trees will connect to a vertex in the other tree,
so the trees $\mathcal{T}$ and $\mathcal{T'}$ cannot both be proper. Thus, for $\mathcal{T} \cap \mathcal{T'} \neq \emptyset$,
trivially \eqref{eqlemsec} holds.

For $\mathcal{T} \cap \mathcal{T'} = \emptyset$, we have to take a closer look at the probabilities involved.
All three probabilities in the lemma are a product over all edges of the probability that either the edge does not connect to any of
the vertices in the tree(s) or the probability that the edge makes a prescribed connection in (one of) the tree(s).
Let $\mathcal{E}_{j,s}(A)$ denote the event that the $j$-th edge of vertex $s$ connects to a vertex in $A$,
with $\mathcal{E}_{j,s}(i)=\mathcal{E}_{j,s}(\{i\})$. Let $\overline{\mathcal{E}}_{j,s}(A)$ be the complement of $\mathcal{E}_{j,s}(A)$.
We have that
\begin{equation}
\mathbb{P}(\mathcal{E}_{j,s}(A)) = \sum_{i\in A} \mathbb{P}(\mathcal{E}_{j,s}(i)),
\end{equation}
because the events on the right-hand side are disjunct. These probabilities are given by the growth rules
\eqref{growthrulePA}--\eqref{growthrulePAc}.

Suppose that the $j$-th edge, $1\leq j\leq m$, of a vertex $t_0$ should not connect to a vertex in $\mathcal{T} \cup \mathcal{T'}$. Then in $\mathbb{P}\left( F_{\mathcal{T}} \cap F_{\mathcal{T'}} \right)$, there will be a factor
\begin{equation}
\mathbb{P}\left( \overline{\mathcal{E}}_{j,t_0}(\mathcal{T}\cup\mathcal{T'})\right)
=1-\mathbb{P}\left(\mathcal{E}_{j,t_0}(\mathcal{T}\cup\mathcal{T'})\right)
=1-\sum_{i\in\mathcal{T} \cup \mathcal{T'}}\mathbb{P}\left(\mathcal{E}_{j,t_0}(i)\right).\label{eqp12}
\end{equation}
In $\mathbb{P}\left( F_{\mathcal{T}}\right)\mathbb{P}\left( F_{\mathcal{T'}}\right)$, there will be a factor
\begin{equation}
\left(1-\sum_{i\in\mathcal{T} }\mathbb{P}\left(\mathcal{E}_{j,t_0}(i)\right)\right)\left(1-\sum_{i\in\mathcal{T'} }\mathbb{P}\left(\mathcal{E}_{j,t_0}(i)\right)\right).\label{eqp1p2}
\end{equation}
It is easy to see that $1-x-y\leq(1-x)(1-y)$ for $x,y\geq0$, so \eqref{eqp1p2} is at least as big as \eqref{eqp12}.

When the $j$-th edge, $1\leq j\leq m$, of a vertex $t_0$, $t+1\leq t_0 \leq 2t$, should connect to a vertex $h \in \mathcal{T}$, then in $\mathbb{P}\left( F_{\mathcal{T}}\cap F_{\mathcal{T'}} \right)$ there will only be a factor
\begin{equation}\label{eqp122}
\mathbb{P}\left(\mathcal{E}_{j,t_0}(h)\right),
\end{equation}
since it will then automatically not connect to a vertex in $\mathcal{T'}$. In $\mathbb{P}\left( F_{\mathcal{T}}\right)\mathbb{P}\left( F_{\mathcal{T'}}\right)$, however, there will be a factor
\begin{equation}
\label{eq122}
\mathbb{P}\left(\mathcal{E}_{j,t_0}(h)\right)\left(1-\sum_{i\in\mathcal{T'}}\mathbb{P}\left(\mathcal{E}_{j,t_0}(i)\right)\right).
\end{equation}
When we multiply \eqref{eq122} by
$\left(1-\sum_{i\in\mathcal{T'}}\mathbb{P}\left(\mathcal{E}_{j,t_0}(i)\right)\right)^{-1}$
we obtain precisely \eqref{eqp122}. By symmetry, the same holds when an edge should connect to a vertex in $\mathcal{T'}$. Since the degree of the vertices in the trees is at most $m+1$, the edges of interest are added after time $t$ and there are less than $m^{k+1}$ vertices in the tree, we have that
\begin{equation}
\left(1-\sum_{i\in\mathcal{T'}}\mathbb{P}\left(\mathcal{E}_{j,t_0}(i)\right)\right)^{-1}\leq \left(1-\frac{m_\delta m^{k+1}}{t}\right)^{-1}.
\end{equation}
Since there are less than $m^{k+1}$ edges in both $\mathcal{T}$ and $\mathcal{T'}$, for $\mathcal{T} \cap \mathcal{T'} = \emptyset$,
\begin{align}
\frac{\mathbb{P}\left( F_{\mathcal{T}} \cap F_{\mathcal{T'}} \right)}
{\mathbb{P}\left(F_{\mathcal{T}}\right)\mathbb{P}\left(F_{\mathcal{T'}}\right)}
&\leq \prod_{h\in \mathcal{T}}
\left(1-\sum_{h\in\mathcal{T'}}\mathbb{P}\left(\mathcal{E}_{j,t_0}(h)\right)\right)^{-1}
\prod_{k\in \mathcal{T'}}
\left(1-\sum_{k\in\mathcal{T'}}\mathbb{P}\left(\mathcal{E}_{j,t_0}(k)\right)\right)^{-1}\nonumber\\
&\leq \left(1-\frac{m_\delta m^{k+1}}{t}\right)^{-2m^{k+1}}=
\left(1+\frac{m_\delta m^{k+1}}{t-m_\delta m^{k+1}}\right)^{2m^{k+1}}\nonumber\\
&\leq \left(1+\frac{m_\delta m \log t}{t-m_\delta m \log t}\right)^{2m \log t}
\leq \left(1+\frac{2m_\delta m \log t}{t}\right)^{2m \log t}.
\end{align}
\end{proof}
We can now use the lemma above to give an upper bound on the variance of $Z_{m,\delta}^{\sss(k)}(2t)$ in terms of the expectation of $Z_{m,\delta}^{\sss(k)}(2t)$.
\begin{proof}[Proof of Proposition~\ref{propsecond}]
Let $c_{m,\delta}=8 m_\delta  m^2$. Then, using Lemma~\ref{lempminpp}, we have that
\begin{align}
\Var\left(Z_{m,\delta}^{\sss(k)}(2t)\right) &= \sum_{\substack{\mathcal{T,T'}\in\mathfrak{T}_m^k(2t) \\ \mathcal{T}\neq{\mathcal{T'}}}}\left(\mathbb{P}\left( F_{\mathcal{T}} \cap F_{\mathcal{T'}} \right) -\mathbb{P}\left(F_{\mathcal{T}}\right)\mathbb{P}\left(F_{\mathcal{T'}}\right)\right)+ \sum_{\mathcal{T}\in\mathfrak{T}_m^k(2t)}\mathbb{P}\left(F_{\mathcal{T}}\right)\left(1-\mathbb{P}\left(F_{\mathcal{T}}\right)\right)\nonumber\\
&\leq \sum_{\substack{\mathcal{T,T'}\in\mathfrak{T}_m^k(2t) \\ \mathcal{T}\neq{\mathcal{T'}}}}\left(\left(1+\frac{2m_\delta m \log t}{t}\right)^{2m \log t} -1\right)\mathbb{P}\left(F_{\mathcal{T}}\right)\mathbb{P}\left(F_{\mathcal{T'}}\right)+ \sum_{\mathcal{T}\in\mathfrak{T}_m^k(2t)}\mathbb{P}\left(F_{\mathcal{T}}\right).\label{eqlema8}
\end{align}
Since
\begin{equation}
\left(1+\frac{2m_\delta m \log t}{t}\right)^{2m \log t} -1 \leq e^{\frac{c_{m,\delta}}{2}\frac{(\log t)^2}{t}}-1 \leq c_{m,\delta} \frac{(\log t)^2}{t},
\end{equation}
we have that \eqref{eqlema8} is at most
\begin{align}
c_{m,\delta} \frac{(\log t)^2}{t} \sum_{\substack{\mathcal{T,T'}\in\mathfrak{T}_m^k(2t) \\ \mathcal{T}\neq{\mathcal{T'}}}}&\mathbb{P}\left(F_{\mathcal{T}}\right)\mathbb{P}\left(F_{\mathcal{T'}}\right) + \mathbb{E}\left[Z_{m,\delta}^{\sss(k)}(2t)\right]\nonumber\\
&\leq c_{m,\delta} \frac{(\log t)^2}{t} \sum_{\mathcal{T,T'}\in\mathfrak{T}_m^k(2t)} \mathbb{P}\left(F_{\mathcal{T}}\right)\mathbb{P}\left(F_{\mathcal{T'}}\right) + \mathbb{E}\left[Z_{m,\delta}^{\sss(k)}(2t)\right]\nonumber\\
&= c_{m,\delta} \frac{(\log t)^2}{t} \mathbb{E}\left[Z_{m,\delta}^{\sss(k)}(2t)\right]^2 + \mathbb{E}\left[Z_{m,\delta}^{\sss(k)}(2t)\right].
\end{align}
\end{proof}

\appendix

\section{Appendix}

\subsection{The tails of the degree sequence}
\begin{lemma}[The total degree of high degree vertices]
\label{totdegree}
Fix $m\geq 1$ and $\delta>-m$.
Assume that  $l_t\to \infty$, as $t\to \infty$ and that
$l_t\leq u_1=t^{\frac{1}{2(\tau-1)}} (\log{t})^{-\frac{1}{2}}$.
Then there exists a constant $B>0$ such that with
probability exceeding $1-o(t^{-1})$,
    \eq
    \label{lowerbddegrees}
    \sum_{i: D_i(t)\geq l_t} D_i(t)\geq B t l_t^{2-\tau}.
    \en
\ch{Moreover, if $N_{\sss \geq l_t}(t)=\#\{i\leq t: D_i(t)\geq l_t\}$ is the number of vertices with degree
at least $l_t$, then, \whp,
    \eq
    N_{\sss \geq l_t}(t) \ch{\geq \sqrt{t}}.
    \en
}
\end{lemma}

\proof We note that
    \eq
    \sum_{i:D_i(t)\geq l_t} D_i(t) \geq  l_t N_{\sss \geq l_t}(t)\ch{.}
    \en

In \cite{DeiEskHofHoo06}, detailed asymptotics for $N_{\sss \geq
l_t}(t)$ were proved for model (c) that we will survey now. These asymptotics play
a key role throughout the proof.

Firstly, it is shown that there exists a $B_1$ such that uniformly
for all $l_t$,
    \eq
    \prob\Big(|N_{\sss \geq l_t}(t)-\expec[N_{\sss \geq l_t}(t)]|
    \geq B_1\sqrt{t\log{t}}\Big)=o(t^{-1}).
    \en
This proves a {\it concentration bound} on the number of vertices
with at least a given degree. The proof of this result follows the
argument in \cite{BolRioSpeTus01}, and holds for any of the models (a)--(c).

Secondly, with
    \eq
    N_{l_t}(t)=\#\{i\leq t: D_i(t)= l_t\},
    \en
the total number of vertices of degree equal to $l_t$, and
with $p_{l_t}$ defined by
    \eq
    \lbeq{pkm2def}
    p_{l_t}=\frac{\ch{(2+\delta/m)}\Gamma(l_t+\delta)\Gamma(m+\delta+\ch{2+\delta/m})}
    {\Gamma(m+\delta)\Gamma(l_t+1+\delta+\ch{2+\delta/m})},\quad l_t\ge m,
    \en
so that $p_k \sim k^{-\tau}$ with $\tau=3+\delta/m$,
there exists a constant $B_2$
such that
    \eq
    \lbeq{sharpBdPl}
    \sup_{l\ge 1} |\expec[N_{l_t}(t)]-tp_{l_t}|\leq B_2.
    \en
For model (c), this is shown in \cite{DeiEskHofHoo06}, for model (a)
this is shown in \cite[Chapter 8]{Hofs08}. This latter proof can easily be adapted
to deal with model (b) as well. In rather generality, results of this
kind (with the sharp bound in \refeq{sharpBdPl}) are proved in \cite{HagWiu06}.

Therefore, we obtain that, with probability exceeding $1-o(t^{-1})$,
    \eqalign
    \label{Pgeqest}
    N_{\sss \geq l_t}(t)
    &\geq \expec[N_{\sss \geq l_t}(t)] -B_1\sqrt{t\log{t}}
    \geq \expec[N_{\sss \geq l_t}(t)]-\expec[N_{\sss \geq 2l_t}(t)]-B_1\sqrt{t\log{t}}\nn\\
    &\geq \sum_{l=l_t}^{2l_t-1} [t\ch{p_{l}}-B_2]-B_1\sqrt{t\log{t}}
    \geq B_3 t l_t^{1-\tau}-B_2l_t-B_1\sqrt{t\log{t}},
    \enalign
for some $B_3 > 0$.
We now wish to pick $l_t$ such that $t l_t^{1-\tau}$ is the dominating term in the
right-hand side of \eqref{Pgeqest}, i.e., $l_t/t^{1/\tau}\to 0$ and
$\sqrt{t\log{t}}/t l_t^{1-\tau} \to 0$, as $t\to \infty$.
Note that $\frac 1 \tau \geq \frac{1}{2(\tau-1)}$ for all $\tau> 2$, so \ch{for} $u_1$ as in the statement of the lemma and $l_t\leq u_1$\ch{,} we find that \eqref{lowerbddegrees} holds with probability exceeding
$1-o(t^{-1})$ and that, \whp, $N_{\sss \geq l_t}(t) \ch{\geq \sqrt{t}}$.
\qed \vskip0.5cm
\subsection{The diameter of the multinomial graph}
\label{multgraph}

\begin{lemma}[Diameter multinomial graph]
\label{lemmamulti}
Let $H_{n_t}$ be the multinomial graph with parameters defined in
\eqref{multinomial}. Then, \whp, the diameter of
$H_{n_t}$ is bounded from above by the
diameter of the uniform Erd\H{o}s-R\'enyi graph $G(n_t,m_t)$, where the
number $m_t$ of edges is equal to
\eq
m_t=\frac12 e_t\Big(1-(1-q_t)^t\Big).
\en
\end{lemma}

\proof
Observe that by definition of the multinomial graph, and with $e_t=n_t(n_t-1)/2$,
    \eq
    \label{rv}
    M_{n_t}=\sum_{i=1}^{e_t} \ch{I[}\sum_{j=1}^t N_{j,i}>0\ch{]}.
    \en
We only have to show that, \whp, the random number of edges
$M_{n_t}$ dominates the deterministic number $m_t$. This can be deduced
from Chebychev's inequality as follows.

From a \ch{straightforward} calculation,
    \eq
    \expec[M_{n_t}]=e_t(1-(1-q_t)^t)=2m_t,
    \en
and
    \eqa
    \Var(M_{\ch{n_t}})&=&
    e_t^2
    \left(
    (1-2q_t)^t-(1-q_t)^{2t}
    \right)
    \ch{ +
    e_t\left(
    (1-q_t)^t-(1-2q_t)^{t}
    \right)}.
    \ena
\ch{The first term is negative, and the second term can be bounded by $e_t(1-(1-q_t)^t)$, so that
    \eq
    \Var(M_{n_t})\leq e_t(1-(1-q_t)^t) = \expec[M_{n_t}] = 2m_t,
    \en
}
so that the variance is of the same order as the first moment.
Applying the Chebychev inequality yields
    \begin{eqnarray}
    \prob(M_{n_t}< m_t)&\le&
    \prob(|M_{n_t}-\expec[M_{n_t}]|>m_t)\leq \frac{\Var(M_{n_t})}{m_t^2}
    \ch{\to 0}.
    \end{eqnarray}
\qed


\subsection{Proof of Lemma \ref{lem-earlypartofcore}} \label{sec-prfslems}

We investigate the problem for model (a) first, the adaptation
of the proof for model (b) is rather straightforward and will be omitted.
The proof for model (c) is slightly more involved
and is treated immediately after the proof for model (a).

We first note that, for models (a) and (b),
the model for general $m\geq 1$ is obtained
from the model for $m=1$ by taking $\delta'=\delta/m$ and
identifying groups of $m$ vertices. For $m=1$ and $\delta>-1$, we shall show by induction on $j$, that for model (a)
    \eqn{
    \lbeq{dibdj=0}
    \prob(D_i(t)=j)\leq C_j\frac{\Gamma(t)\Gamma(i+\Delta)}
    {\Gamma(t+\Delta)\Gamma(i)},
    }
for all $t\geq i$ and $j\geq m$, with $\Delta=(1+\delta)/(2+\delta)\in(0,1)$
and where $C_j$ will be determined in the course of the proof.
Clearly, for every $t\geq i$, for model (a),
    \eqn{
    \lbeq{prod-bd}
    \prob(D_i(t)=1)=\prod_{s=i+1}^{t} \Big(1-\frac{1+\delta}{(2+\delta)(s-1)+(1+\delta)}\Big)
    =\frac{\Gamma(t)\Gamma(i+\Delta)}
    {\Gamma(t+\Delta)\Gamma(i)},
    }
which initializes the induction hypothesis with $C_1=1$.

To advance the induction, we let $s\leq t$ be the last time at
which a vertex is added to $i$. Then we have that
    \eqn{
    \lbeq{bdIndaa}
    \prob(D_i(t)=j)=\sum_{s=i+j-1}^t \prob\big(D_i(s-1)=j-1\big)\frac{j-1+\delta}{(2+\delta)(s-1)+1+\delta}
    \prob\big(D_i(t)=j|D_i(s)=j\big).
    }
By the induction hypothesis, we have that
    \eqn{
    \lbeq{bdInda}
    \prob\big(D_i(s-1)=j-1\big)\leq C_{j-1}\frac{\Gamma(s-1)\Gamma(i+\Delta)}
    {\Gamma(s-1+\Delta)\Gamma(i)}.
    }
Moreover, analogously to \refeq{prod-bd}, we have that
    \eqan{
    \lbeq{bdIndb}
    \prob(D_i(t)=j|D_i(s)=j)&=\prod_{q=s+1}^{t} \Big(1-\frac{j+\delta}{(2+\delta)(q-1)+(1+\delta)}\Big)\\
    &
    =\frac{\Gamma(t-\frac{j-1}{2+\delta})\Gamma(s+\Delta)}
    {\Gamma(t+\Delta)\Gamma(s-\frac{j-1}{2+\delta})}.\nn
    }
Combining \refeq{bdIndaa}, \refeq{bdInda} and \refeq{bdIndb}, we arrive at
    \eqan{
    \prob(D_i(t)=j)&\leq C_{j-1}\sum_{s=i+j-1}^t \frac{\Gamma(s-1)\Gamma(i+\Delta)}
    {\Gamma(s-1+\Delta)\Gamma(i)}
    \frac{j-1+\delta}{(2+\delta)(s-1)+(1+\delta)}
    \frac{\Gamma(t-\frac{j-1}{2+\delta})\Gamma(s+\Delta)}
    {\Gamma(t+\Delta)\Gamma(s-\frac{j-1}{2+\delta})}\nonumber\\
    &=C_{j-1}\frac{j-1+\delta}{2+\delta} \frac{\Gamma(i+\Delta)}
    {\Gamma(i)}\frac{\Gamma(t-\frac{j-1}{2+\delta})}{\Gamma(t+\Delta)}
    \sum_{s=i+j-1}^t  \frac{\Gamma(s-1)}
    {\Gamma(s-\frac{j-1}{2+\delta})}.
    }
We note that, whenever $l+b, l+1+a>0$ and $a-b+1>0$,
    \eqn{
    \lbeq{Gammasum}
    \sum_{s=l}^t  \frac{\Gamma(s+a)}
    {\Gamma(s+b)}=\frac{1}{a-b+1} \Big[\frac{\Gamma(t+1+a)}{\Gamma(t+b)}-\frac{\Gamma(l+a)}{\Gamma(l-1+b)}\Big]
    \leq \frac{1}{a-b+1}\frac{\Gamma(t+1+a)}{\Gamma(t+b)}.
    }
Application of \refeq{Gammasum} for $a=-1, b=-\frac{j-1}{2+\delta}, l=i+j-1$,
so that $a-b+1=\frac{j-1}{2+\delta}>0$ when $j>1$,
leads to
    \eqan{
    \lbeq{IndAdv}
    \prob(D_i(t)=j)
    &\leq C_{j-1}\frac{j-1+\delta}{j-1}\frac{\Gamma(i+\Delta)}
    {\Gamma(i)}\frac{\Gamma(t)}
    {\Gamma(t+\Delta)}.
    }
Equation \refeq{IndAdv} advances the induction when we define
\ch{
    \eqn{
    \lbeq{Cj-def}
    C_j=\frac{\Gamma(j+\delta)}{\Gamma(j)\Gamma(1+\delta)},
    }
so that}
    \eqn{
    C_j=\frac{j-1+\delta}{j-1}C_{j-1}.
    }
For $m>1$, inequality \refeq{dibdj=0} for model (a) generalizes to
    \eqn{
    \prob(D_i(t)=j)\leq C_j\frac{\Gamma(t)\Gamma(i+\frac{1+\delta'}{2+\delta'})}
    {\Gamma(t+\frac{1+\delta'}{2+\delta'})\Gamma(i)}
    =C_j\frac{\Gamma(t)\Gamma(i+\frac{m+\delta}{2m+\delta})}
    {\Gamma(t+\frac{m+\delta}{2m+\delta})\Gamma(i)}.
    }
This completes the investigation
of $\prob(D_i(t)=j)$ for model (a).
In an identical fashion, for model (b), we obtain for $m=1$
    \eqn{
    \lbeq{dibdj=0(b)}
    \prob(D_i(t)=j)\leq C_j\frac{\Gamma(t-\Delta)\Gamma(i)}
    {\Gamma(t)\Gamma(i-\Delta)},
    }
where again $C_1=1$ and $C_j$ satisfies \refeq{Cj-def}. This generalizes to
    \eqn{
    \prob(D_i(t)=j)\leq C_j\frac{\Gamma(t-\frac{m+\delta}{2m+\delta})\Gamma(i)}
    {\Gamma(t)\Gamma(i-\frac{m+\delta}{2m+\delta})}.
    }
We omit further details for
model (b).

For models (a) and (b) we can generalize the inequality for $m=1$
to $m>1$. Unfortunately this fails for model (c), and we first adapt the argument.
Recall that $D_i(t)$ is the degree of vertex $i$ at time $t$.
We shall define $E_i(t)$ such that $E_i(t)\leq D_i(t)$ and
$E_i(t)$ grows by at most one at each time step. The definition of
$E_i(t)$ is recursive. We let $E_i(i)=D_i(i)=m$, and, assuming we have shown
that $D_i(t)=E_i(t)+R_i(t)$, where $R_i(t)\geq 0$, we proceed at time $t+1$
as follows. We can increase $E_i(t)$ only when the {\it first} edge of
vertex $t+1$ attaches to vertex $i$, and we do this with probability
$\frac{E_i(t)+\delta}{(2m+\delta)t}$. With probability $\frac{R_i(t)}{(2m+\delta)t}$,
we keep $E_i(t+1)=E_i(t)$ and we increase $R_i(t)$ by one. For the other $m-1$
edges, we increase $R_i(t)$ by one with probability $\frac{D_i(t)+\delta}{(2m+\delta)t}$.
Then we clearly have that $E_i(t+1)\leq D_i(t+1)$ if $E_i(t)\leq D_i(t)$, since
the difference between $D_i(t)$ and $E_i(t)$ equals $R_i(t)$, which is
monotonically increasing. Moreover, we have that $E_i(t+1)$ equals
$E_i(t)$ or $E_i(t)+1$, and the latter occurs with conditional probability
    \eqn{
    \prob(E_i(t+1)=j|E_i(t)=j-1)=\frac{j-1+\delta}{(2m+\delta)t}.
    }
We now adapt the above argument for model (a) to the random variable
$E_i(t)$. Indeed, we
now use as an induction hypothesis that
    \eqn{
    \lbeq{dibdj=0(c)}
    \prob(E_i(t)=j)\leq C_j\frac{\Gamma(t-\frac{m+\delta}{2m+\delta})\Gamma(i)}
    {\Gamma(t)\Gamma(i-\frac{m+\delta}{2m+\delta})},
    }
where $C_m=1$ and, for $j>m$,
    \eqn{
    \lbeq{Cj(c)}
    C_j=\frac{j-1+\delta}{j-m}C_{j-1}.
    }
The verification of \refeq{dibdj=0(c)} is a straightforward adaptation of the
one of \refeq{dibdj=0}.

We summarize the bounds in models (a)--(c):
for all $m\geq 1$, and $i\in [t], j\geq m$,
    \eqn{
    \lbeq{summdjbd}
    \prob(\ch{E_i}(t)=j)\leq C_j\frac{\Gamma(t-a_1)\Gamma(i+a_2)}
    {\Gamma(t+a_2)\Gamma(i-a_1)},
    }
\ch{where $E_i(t)=D_i(t)$ in models (a) and (b) and} where $a_1=0$ for model (a), while $a_1=\frac{m+\delta}{2m+\delta}$
for models (b)--(c), while $a_2=\frac{m+\delta}{2m+\delta}$ for model (a),
while $a_2=0$ for models (b)--(c), and, for all models,
$C_j\leq j^{p-1}$ for some $p\geq 1$.

Consequently, we obtain
    \eqn{
    \lbeq{consdibdj=0}
    \prob(D_i(t)\leq j) \leq j^{p} \frac{\Gamma(t-a_1)\Gamma(i+a_2)}
    {\Gamma(t+a_2)\Gamma(i-a_1)}.
    }
We finally use \refeq{consdibdj=0} to complete the proof of Lemma \ref{lem-earlypartofcore}.
Take $0<b<\frac{a_1+a_2}{a_1+a_2+1}=\frac{m+\delta}{3m+2\delta}$, then, by Boole's inequality,
    \eqan{
    \prob(\exists i\leq t^b: D_i(t)\leq (\log{t})^{\sigma})
    &\leq \sum_{i=1}^{t^b} \prob(D_i(t)\leq
    (\log t)^{\sigma})
    \leq (\log{t})^{\sigma p} \frac{\Gamma(t-a_1)}
    {\Gamma(t+a_2)}\sum_{i=1}^{t^b} \frac{\Gamma(i+a_2)}
    {\Gamma(i-a_1)}\nn\\
    &\leq (\log{t})^{\sigma p}(a_1+a_2+1)^{-1}\frac{\Gamma(t-a_1)}
    {\Gamma(t+a_2)}\frac{\Gamma(t^b+a_2+1)}
    {\Gamma(t^b-a_1)}
    =o(1).
    }
This completes the proof of Lemma \ref{lem-earlypartofcore}.
\qed

\subsection{Late vertices have small degree}
Recall the definition of the core $\Core$ in \eqref{defcore}, where we take $\sigma > 1$.
In the following theorem we will prove that, for models (a)--(c), all vertices with large degree will be early vertices.
We need this result to prove Theorem~\ref{thmlowerb}.
\begin{proposition}[Late vertices have small degree]
\label{latevertsmalldeg}
Fix $m\geq 2, \delta >-m$ and $\sigma>1$. Then,
${\rm Core}_{2t}\subseteq [t]$ \whp.
\end{proposition}

\begin{proof}
Note that
\begin{align}
\mathbb{P}\Big({\rm Core}_{2t}\subseteq[t]\Big) \geq1-\sum_{i=t+1}^{2t} \mathbb{P}\Big(D_i(2t)\geq(\log 2t)^\sigma\Big)
&\geq 1-\sum_{i=t+1}^{2t}\mathbb{P}\Big(D_t(2t)\geq(\log 2t)^\sigma\Big)\nonumber\\
&= 1-t\mathbb{P}\Big(D_t(2t)\geq(\log 2t)^\sigma\Big),
\end{align}
because vertex $t$ is more likely to have a large degree than vertices added after time $t$. In Lemma~\ref{lempolya} we will show that $\mathbb{P}\big(D_t(2t)\geq(\log 2t)^\sigma\big)=o\left(\frac{1}{t}\right)$, so that
$\mathbb{P}\Big({\rm Core}_{2t}\subseteq[t]\Big) \geq 1 - o(1).$
\end{proof}

\begin{lemma}[Tails of degree distribution]
\label{lempolya} Fix $m\geq 2, \delta >-m$ and $\sigma>1$. Then,
\begin{equation}\label{eqlempolya}
\mathbb{P}\Big(D_{t}(2t)\geq(\log 2t)^\sigma\Big)=o\left(1/t\right).
\end{equation}
\end{lemma}
\begin{proof}
\ch{We investigate the problem for models (a) and (b) first, the adaptation
for model (c) will be discussed later.}
As noted in Section~\ref{sec-intro}, \ch{for models (a) and (b),} $G_{m,\delta}(2t)$ can be constructed from $G_{1,\delta'}(2mt)$, with $\delta'=\delta/m$. We will include the superscripts to avoid confusion. Thus identify, for $i\in [2t]$, vertices $((i-1)m+1)^{\sss(1)},\ldots,(im)^{\sss(1)}$ in $G_{1,\delta'}(2mt)$ with vertex $i^{\sss(m)}$ in $G_{m,\delta}(2t)$. So \eqref{eqlempolya} is equivalent to
\begin{equation}
\label{eqvssmall}
\mathbb{P}\Big(D_{((t-1)m+1)^{\sss(1)}}(2mt)+\ldots+D_{(tm)^{\sss(1)}}(2mt)\geq(\log 2t)^\sigma\Big)=o\left(1/t\right).
\end{equation}

We will now color the vertices and edges in the following way. Color the vertices $1^{\sss(1)}, \ldots, ((t-1)m)^{\sss(1)}$ and all edges between these vertices blue and color the vertices $((t-1)m+1)^{\sss(1)},\ldots,(tm)^{\sss(1)}$ and the $m$ edges that are attached to them at time $mt$ red. When a vertex that was added after time $mt$ connects to a blue (red) vertex, also color that vertex and its edge blue (red). Color vertices with a self-loop and its edge blue. Then, at time $2mt$, the total degree of vertices $((t-1)m+1)^{\sss(1)},\ldots,(tm)^{\sss(1)}$ is at most equal to the number of red edges plus $m$, because no blue edges are connected to these red vertices, and all red edges are connected with at most one endpoint to these vertices. The only exception are the first $m$ red edges, which might connect with both endpoints to these vertices, hence we have to add $m$ to the number of red edges. Thus,
    \begin{equation}
    \mathbb{P}\Big(D_{((t-1)m+1)^{\sss(1)}}(2mt)+\ldots+D_{(tm)^{\sss(1)}}(2mt)\geq(\log 2t)^\sigma\Big) \leq \mathbb{P}\Big(\text{\#\{red edges\} } + m\geq(\log 2t)^\sigma\Big).
    \end{equation}
Since we will bound the right-hand side of the formula above, it is allowed to increase the probability of attaching to a red vertex, or, equivalently, to decrease the probability of attaching to a blue vertex. It is also allowed to increase the total degree of the red vertices, or to decrease the total degree of the blue vertices. All this will only increase the probability of the number of red edges being large.

Therefore, we are allowed to assume that the first $m$ red edges are all self-loops. Further, we will not allow for self-loops after time $t$, which will increase the probability of attaching to a red vertex in models (a) and (b), in model (c) nothing changes. When we consider model (c), we see that the degrees should only be updated after each $m$-th vertex has been added. For $j \geq mt$, no more than $m$ edges and vertices can be added before updating the degrees, so
    \begin{align}
    \mathbb{P}\Big((j+1)^{\sss(1)} \text{ connects to a red vertex}\big| G_{1,\delta'}^{(c)}(j)\Big) &= \frac{\sum_{v^{\sss(1)} \text{ red}}\left( D_{v^{\sss(1)}}(m\lfloor j/m\rfloor)+\delta'\right)}{m\lfloor j/m \rfloor (2+\delta')} \nonumber\\
    &\leq \frac{\sum_{v^{\sss(1)} \text{ red}}\left( D_{v^{\sss(1)}}(j)+\delta'\right)}{j (2+\delta')-m(2+\delta')}.
    \end{align}
Thus, we are allowed to update the degrees after adding each vertex, but then we have to lower the total weight that blue vertices and edges contribute to the connecting probabilities by $m(2+\delta')$. The above bound on the connecting probabilities also holds for models (a) and (b).

Since we are only interested in the number of red and blue vertices and edges, the problem reduces to the following P\'olya urn scheme. Let there be an urn with, at time $s$, $S_1(s)$ red balls, corresponding to the total weight that red vertices and edges contribute to the connecting probabilities, and $S_2(s)$ blue balls, corresponding to the lowered total weight that blue vertices and edges contribute to the connecting probabilities. At time $s=0$ we will start with $S_1(0)=m(2+\delta')$ and $S_2(0) = m(t-1)(2+\delta')-m(2+\delta')$. We then successively take one ball proportional to the number of balls of a certain color, and replace it together with another $2+\delta'$ balls of the same color. This corresponds to attaching a new vertex to a vertex of that color.

So $\frac{S_1(mt)}{2+\delta'}$ has the same distribution as the number of red edges at time $2mt$. Consequently,
    \begin{equation}
    \mathbb{P}\Big(D_{t^{\sss(m)}}(2t)\geq(\log 2t)^\sigma\Big) \leq \mathbb{P}\left(\frac{S_1(mt)}{2+\delta'}+m \geq (\log 2t)^\sigma\right).
    \end{equation}
To analyse the probability on the right-hand side, we make use of De Finetti's Theorem \cite{Fel70}. This theorem states that for an infinite sequence of exchangeable random variables $\{X_i\}_{i=1}^{\infty}, X_i\in\{0,1\}$, there exists a random variable $U$ with $\mathbb{P}(U \in [0,1])=1$, such that for all $1\leq k \leq n$,
    \begin{equation}
    \mathbb{P}\Big(X_1 = \ldots = X_k = 1, X_{k+1} =0, \ldots, X_n=0\Big) = \mathbb{E}\Big[U^k(1-U)^{n-k} \Big].
    \end{equation}
The random variable $U$ can be computed explicitly. Note that this implies that
    \begin{equation}
    \mathbb{P}\left( \sum_{i=1}^{n} X_i = k \right) = \mathbb{E}\Big[ \mathbb{P}\Big({\rm BIN}(n, U) = k \Big| U \Big)\Big].
    \end{equation}

Let $X_i$ denote the indicator that the $i$-th ball drawn in the P\'olya urn scheme described above is red.
As shown in \cite[Section~11.1]{Hofs08},
$\{X_i\}_{i=1}^{\infty}$ is an infinite exchangeable sequence. Note that
    \begin{equation}
    S_1(s) = (2+\delta')m + (2+\delta')\sum_{i=1}^{s} X_i.
    \end{equation}
Hence,
    \begin{equation}
    \mathbb{P}\left(\frac{S_1(mt)}{2+\delta'}+m \geq (\log 2t)^\sigma\right)=\mathbb{E}\Big[\psi(U)\Big],\label{eqebin}
    \end{equation}
where $0 \leq \psi(u) = \mathbb{P}\Big({\rm BIN}(mt,u) \geq (\log 2t)^\sigma-2m\Big) \leq 1$.

Now observe from \cite{JanLucRuc00} that
\begin{equation}\label{eqpsiu}
\psi(u) \leq e^{-(\log 2t)^{\sigma}+2m},
\end{equation}
whenever $u$ is such that $7mtu \leq (\log 2t)^{\sigma}-2m$. We define $g(t)=((\log 2t)^\sigma-2m)/(7(mt))$. Since,
\begin{align}
\mathbb{E}\Big[\psi(U)\Big] &= \mathbb{E}\Big[\psi(U)\Big|U \leq g(t)\Big]\mathbb{P}\Big(U \leq g(t)\Big)+\mathbb{E}\Big[\psi(U)\Big|U > g(t)\Big]\mathbb{P}\Big(U > g(t)\Big) \nonumber\\
&\leq \psi(g(t))+\mathbb{P}\Big(U > g(t)\Big),
\end{align}
we obtain, according to \eqref{eqpsiu},
\begin{equation}
\mathbb{P}\left(\frac{S_1(mt)}{2+\delta'}+m \geq (\log 2t)^\sigma\right) \leq e^{-(\log 2t)^{\sigma}+2m} + \mathbb{P}\Big(U > g(t)\Big) = o\left(\frac{1}{t}\right) + \mathbb{P}\Big(U > g(t)\Big).
\end{equation}
It remains to show that also $\mathbb{P}(U> g(t))=o\left(\frac{1}{t}\right)$. It turns out that $U$ has a Beta-distribution with parameters $\alpha=m$ and $\beta=m(t-2)$ (\cite{Hofs08}), so $\alpha, \beta>1$. Thus we have that the probability density function of $U$ is unimodular, with its turning point at $t=\frac{\alpha-1}{\alpha+\beta-2}$ (\cite{WadBry60}). It is easy to verify that $g(t)\geq \frac{\alpha-1}{\alpha+\beta-2}$, for $t$ sufficiently large, so that
\begin{align}
\mathbb{P}(U>g(t)) &\leq \left(1-g(t)\right) \frac{\Gamma(\alpha+\beta)}{\Gamma(\alpha)\Gamma(\beta)} \left(g(t)\right)^{\alpha-1} \left(1-g(t)\right)^{\beta-1}
\leq \frac{\Gamma(\alpha+\beta)}{\Gamma(\alpha)\Gamma(\beta)} \left(1-g(t)\right)^{\beta}.\label{eqpuggt}
\end{align}
Using Stirling's formula (see e.g., \cite{AbrSte64}), one can show that  there exists a constant $C>0$, such that \eqref{eqpuggt} is at most
\begin{align}
C\frac{\beta^{\alpha}}{\Gamma(\alpha)}\left(1-g(t)\right)^{\beta} &\leq C(mt)^{m}\left(1-\frac{(\log 2t)^{\sigma}}{8m(t-2)}\right)^{m(t-2)}
\leq C(mt)^{m}e^{-(\log 2t)^{\sigma}/8}
=o\left(1/t\right),
\end{align}
because $\sigma > 1$.

Note that we in fact proved that $\mathbb{P}(D_{t^{\sss(m)}}(2t)\geq(\log 2t)^\sigma)=o\left(t^{-\gamma}\right)$, for any constant $\gamma$.
\end{proof}

\paragraph*{Acknowledgements.}
The work of RvdH and SD is supported in part by Netherlands Organisation for
Scientific Research (NWO). We thank Mia Deijfen and Henri van den
Esker for many useful conversations throughout the project\ch{, and the referee for many remarks improving the presentation of the paper}.

\end{document}